\documentclass[12pt,a4paper]{article}  
\usepackage{amsmath,amsfonts,amsthm,amssymb}
\usepackage[all]{xy}
                             
\title{Deformations of Type D Kleinian Singularities}  
\author{Paul Boddington}     
\date{December 2006}

\begin{document}             

\maketitle                 

\newtheorem{lemma}{Lemma}[section]
\newtheorem{defi}[lemma]{Definition}
\newtheorem{theo}[lemma]{Theorem}
\newtheorem{prop}[lemma]{Proposition}
\newtheorem{cor}[lemma]{Corollary}
\newtheorem{rema}[lemma]{Remark}

\begin{abstract}
For $n\geq 4$ we shall construct a family $D(q)$ of non-commutative deformations of the coordinate algebra of a Kleinian singularity of type $D_n$ depending on a polynomial $q$ of degree $n$. We shall prove that every deformation of a type $D$ Kleinian singularity that is not commutative is isomorphic to some $D(q)$. We shall then consider in type $D$ the family of deformations $\mathcal{O}^{\boldsymbol{\lambda}}$ constructed by Crawley-Boevey and Holland. For each $\mathcal{O}^{\boldsymbol{\lambda}}$ that is not commutative we shall exhibit an explicit isomorphism $D(q)\cong \mathcal{O}^{\boldsymbol{\lambda}}$ for a suitable choice of $q$. This will enable us to prove that every deformation of a Kleinian singularity of type $D_n$ that is not commutative is isomorphic to some $\mathcal{O}^{\boldsymbol{\lambda}}$ and determine when two $\mathcal{O}^{\boldsymbol{\lambda}}$ are isomorphic.
\end{abstract}
\section{Introduction}
Kleinian singularities are in bijective correspondence with the extended Dynkin diagrams $A_{n-1}$ ($n\geq2$), $D_n$ ($n\geq4$), $E_6$, $E_7$ and $E_8$. Deformations of the coordinate algebra of Kleinian singularities are becoming increasingly important in representation theory and mathematical physics. Suppose $\mathfrak{g}$ is a complex simple Lie algebra and $E$ is a subregular nilpotent element. Extend this to an $\mathfrak{sl}_2$ triple $E$, $H$, $F$. Brieskorn \cite{brieskorn} showed that the intersection of the Slodowy slice \cite{slodowy} $E+Z_{\mathfrak{g}}(F)$ with the nilpotent cone of $\mathfrak{g}$ is isomorphic to a Kleinian singularity. In fact one can carry out the same process for an arbitrary nilpotent orbit, and Premet showed in \cite{premet} that all the singularities constructed in this way admit natural deformations. Premet's deformations have been shown to be isomorphic to finite $W$-algebras (See \cite{a}, \cite{SK}, \cite{boer}). In \cite{gordonrumynin} Gordon and Rumynin consider connections between deformations and modular representation theory. 
\par
There is a bijective correspondence, due to McKay \cite{mckay}, between finite non-trivial subgroups $G$ of $SL_2(\mathbb{C})$ and the extended Dynkin diagrams. For such a group, the quotient $\mathbb{C}^2/G$ is a Kleinian singularity and has coordinate algebra $\mathbb{C}[x,y]^{G}$ where $G$ acts on $\mathbb{C}[x,y]$ in the obvious way. In \cite{CBH}, Crawley-Boevey and Holland use related ideas to construct a family $\mathcal{O}^{\boldsymbol{\lambda}}$ of non-commutative deformations of the coordinate algebra of a Kleinian singularity. The construction is as follows. $G$ acts on the algebra of non-commuting polynomials so we can form the skew group algebra $\mathbb{C}\langle x,y\rangle\ast G$. Next, pick $\boldsymbol{\lambda}\in Z(\mathbb{C}G)$. $\mathcal{O}^{\boldsymbol{\lambda}}$ is then defined to be the algebra $e(\mathbb{C}\langle x,y\rangle\ast G/(xy-yx-\boldsymbol{\lambda}))e$ where $e$ is the average of the group elements. $\mathcal{O}^{\boldsymbol{\lambda}}$ is commutative if and only if the trace of $\boldsymbol{\lambda}$ on the regular representation of $G$ is zero. The algebras we shall call $\mathcal{O}^{\boldsymbol{\lambda}}$ have a different description, also due to Crawley-Boevey and Holland. It follows from Theorem 3.4 of \cite{CBH} that the two families can be identified. 
\par
For a Kleinian singularity of type $A_{n-1}$ the coordinate algebra can be defined by the simple equation $u^n=vw$ and a family of deformations $T(s)$ was constructed by Hodges (see \cite{hodges}), using ideas similar to those of Smith in \cite{smith}. These algebras had in fact been studied earlier by Bavula \cite{bav1}, \cite{bav2}. $T(s)$ is defined in terms of generators and relations and depends on a polynomial $s$. Bavula and Jordan \cite{bav} gave necessary and sufficient for two algebras $T(s)$ and $T(s')$ to be isomorphic. It is not too difficult to see that every algebra $T(s)$ is isomorphic to one of the $\mathcal{O}^{\boldsymbol{\lambda}}$ that is not commutative, and vice-versa.
\par
In type $D$ the situation is more complicated. The problem of describing all deformations of (the coordinate algebra of) a Kleinian singularity of type $D_n$ $(n\geq 4)$ in terms of generators and relations, and determining when two are isomorphic, wasn't solved until recently. The complete family was first described by generators and relations in \cite{boddington}, where one can also find the solution to the isomorphism problem for $n\geq 5$. The family was later discovered independently by Levy \cite{levy} who also solved the isomorphism problem for $n\geq 4$. Levy's result is stronger in that his classification is up to isomorphism of algebras, whereas the earlier result was only a classification up to filtered algebras. Although the parameterizations used by the two authors were slightly different, both families depended on a polynomial of degree $n-1$ and a complex number. It is not at all obvious from the presentations of the deformations how the algebras relate to the $\mathcal{O}^{\boldsymbol{\lambda}}$. 
\par
Here, we shall work exclusively in type $D_n$. We shall give a new construction of the family of deformations, which we shall call $D(q)$ where $q$ is a polynomial of degree $n$. The construction will exhibit $D(q)$ as a subalgebra of an algebra similar to those considered by Bavula, Smith and Hodges. For each $\mathcal{O}^{\boldsymbol{\lambda}}$ we shall describe an explicit isomorphism $D(q)\cong\mathcal{O}^{\boldsymbol{\lambda}}$ for a suitable choice of $q$. This will enable us to prove that every non-commutative deformation of a Kleinian singularity of type $D_n$ is isomorphic as a filtered algebra to some $\mathcal{O}^{\boldsymbol{\lambda}}$, and determine all the isomorphisms between different $\mathcal{O}^{\boldsymbol{\lambda}}$. It turns out that the only possible isomorphisms occur from graph automorphisms, scaling $\boldsymbol{\lambda}$, and reflections (see section 8).
\par
In type $D_n$ the coordinate algebra is usually defined by the equation $x^{n-1}+xy^2+z^2=0$, but we will find it more convenient to use the equation $(-1)^n u^{n-1}-uv^2+w^2=0$ which is clearly equivalent under an automorphism of the polynomial algebra to the standard one. The coordinate algebra is graded by putting $u$, $v$, $w$ in degrees $4$, $2n-4$ and $2n-2$ respectively. We shall adopt the convention that $0\in \mathbb{N}$. A final point about notation is that $p'$ shall never mean the derivative of a polynomial $p$.
\par
\emph{Acknowledgements} This project began while I was a PhD student at the University of Warwick and I would like to thank Dmitriy Rumynin for his guidance. I would also like to thank Mark Boddington for the help he has given me in eliminating errors, and Paul Levy for bringing \cite{bav1} and \cite{bav2} to my attention.  

\section{Algebras $\tilde{T}(s)$}      % Produces section heading.  Lower-level
                             % sections are begun with similar 
                             % \subsection and \subsubsection commands.

Let $G$ be an infinite cyclic group generated by an element $a$, and let $\mathbb{C}(h)$ denote the $\mathbb{C}$-algebra of rational functions in an indeterminate $h$. An action of $G$ on $\mathbb{C}(h)$ is given by \[t(h)^{a^k} = t(h+k)\text{ for all }t(h)\in \mathbb{C}(h),k\in\mathbb{Z}.\] We may therefore form the skew group algebra $G\ast \mathbb{C}(h)$. By definition, $G\ast \mathbb{C}(h)$ is the complex vector space $\mathbb{C}G\otimes \mathbb{C}(h)$ with multiplication determined by
\[(a^i\otimes p(h))(a^j \otimes q(h))=a^{i+j} \otimes p(h+j)q(h)\text{ for all }i,j\in \mathbb{Z}, p(h),q(h)\in\mathbb{C}(h).\]
From now on we shall drop the $\otimes$ symbols.
\begin{defi}
Suppose $s(h)\in\mathbb{C}(h)\backslash \{0\}$. Define $\tilde{T}(s)$ to be the algebra $G\ast \mathbb{C}(h)$. Define also $b:=s(h)a^{-1}\in\tilde{T}(s)$.
\end{defi}
\begin{rema}
The algebra $\tilde{T}(s)$ is clearly independent of $s$, but since we shall regard the element $b$ as being part of the structure, the name makes sense.
\end{rema}
\begin{lemma}
In $\tilde{T}(s)$ we have $ba=s(h)$ and $ab=s(h-1)$. Moreover for all $t(h)\in \mathbb{C}(h)$, $k\in\mathbb{Z}$ we have
\[t(h)a^k = a^k t(h+k), \hspace{10mm} t(h)b^k = b^k t(h-k).\]
\end{lemma}
\noindent
Proof. Obviously $ba=s(h)a^{-1}a=s(h)$. Next, the fact that $t(h)a^k = a^k t(h+k)$ follows from the definition of the skew group algebra. It follows that 
\[ab=as(h)a^{-1}=aa^{-1}s(h-1)=s(h-1).\] Next, \[t(h)b=t(h)s(h)a^{-1}=s(h)t(h)a^{-1}=s(h)a^{-1}t(h-1)=bt(h-1),\] so the result for $t(h)b^k$ for nonnegative $k$ follows by induction. Similarly, \[t(h)b^{-1}=t(h)a\frac{1}{s(h)}=at(h+1)\frac{1}{s(h)}=a\frac{1}{s(h)}t(h+1)=b^{-1}t(h+1).\]
The result follows. $\Box$

\begin{rema}
If $s$ is a polynomial of degree $n$, then the subalgebra of $\tilde{T}(s)$ generated by $a$, $b$ and $h$ is denoted $T(s)$ and can be described by the presentation $ha=a(h+1)$, $hb=b(h-1)$, $ab=s(h-1)$, $ba=s(h)$. It follows that if $T(s)$ is given a filtration by putting $a$, $b$ and $h$ in degrees $n$, $n$, $2$ respectively then $T(s)$ is a deformation of a Kleinian singularity of type $A_{n-1}$. 
\end{rema}
\noindent The algebras $T(s)$ were considered by Bavula \cite{bav1}, \cite{bav2} and Hodges \cite{hodges}. 

\section{Some Notation}

In this section we shall introduce some notation and collect together some simple results involving the notation. We shall make heavy use of these ideas in what follows, as the notation will appear in our presentations of the algebras $D(q)$.  

\begin{defi}
Suppose $f(\sqrt{x} )$ is a polynomial in $\sqrt {x}$ where $x$ is a formal indeterminate and $\sqrt{x}$ a formal square root. Then $f(\sqrt{x} )$ can be written uniquely in the form
\[ f(\sqrt{x}) = g(x) + \sqrt{x} h(x),\]
where $g(x)$ and $h(x)$ are polynomials in $x$. Define $\big\{ f(\sqrt{x} )\big\}:=h(x)$.
More generally, suppose $r$ is any element of any $\mathbb{C}$-algebra. Define $\big\{ f(\sqrt{r}) \big\}$ to be the result of substituting $r$ for $x$ in $\big\{ f(\sqrt{x} )\big\}$. There is no assumption here that $r$ has a square root. 
\end{defi}

\begin{lemma}
\label{root2} For any polynomial $f$ we have
\[\big\{ f(\sqrt{x} )\big\} = \frac{f(\sqrt{x})-f(-\sqrt{x})}{2\sqrt{x}}, \hspace{5mm} \big\{ \sqrt{x}f(\sqrt{x} )\big\} = \frac{f(\sqrt{x})+f(-\sqrt{x})}{2}. \]
\end{lemma}

\noindent
Proof. Follows from the equality
\[f(\sqrt{x}) = \frac{f(\sqrt{x})+f(-\sqrt{x})}{2} + \sqrt{x}\frac{f(\sqrt{x})-f(-\sqrt{x})}{2\sqrt{x}}.\Box\]

\noindent
The following lemma is obvious.

\begin{lemma}
$\left\{ \cdot\right\}$ is $\mathbb{C} [x]$-linear in the sense that if $f$, $g$ and $p$ are polynomials, then $\left\{ f( \sqrt{x})+g( \sqrt{x}) \right\}=\left\{f( \sqrt{x})\right\} + \left\{g( \sqrt{x})\right\}$ and $\big\{ p(x)f(\sqrt{x}) \big\} = p(x)\big\{ f(\sqrt{x}) \big\}.$
 $\Box$  
\end{lemma}

\begin{lemma}
\label{bridge}
For any polynomial $f$ we have $\left\{f(\sqrt{x})^2\right\}=2\left\{f(\sqrt{x})\right\}\left\{\sqrt{x}f(\sqrt{x})\right\}$ and $\left\{ \sqrt{x}f(\sqrt{x})^2\right\}=x\left\{f(\sqrt{x})\right\}^2 +\left\{\sqrt{x}f( \sqrt{x})\right\}^2$.
\end{lemma}

\noindent
Proof. Write $f(\sqrt{x}) = g(x) + \sqrt{x} h(x)$. Then $\left\{f(\sqrt{x})^2\right\}=\left\{ (g(x)+\sqrt{x}h(x))^2\right\}=2g(x)h(x)=2\left\{f(\sqrt{x})\right\}\left\{\sqrt{x}f(\sqrt{x})\right\}$. Similarly $\left\{ \sqrt{x}f(\sqrt{x})^2\right\} = \left\{ \sqrt{x}(g(x) +\sqrt{x}h(x))^2\right\}=g(x)^2+xh(x)^2$. $\Box$
\begin{lemma}
\label{iso}
For any $\lambda\in \mathbb{C}$, the map $f( \sqrt{x})\mapsto \left\{(\sqrt{x} - \lambda )f( \sqrt{x})\right\}$ is a vector space isomorphism from the space of polynomials in $\sqrt{x}$ satisfying $f(-\sqrt{x})=f(\sqrt{x}-1)$ to the space $\mathbb{C} [x]$.
\end{lemma}

\noindent
Proof. We first claim that the condition that $f(-\sqrt{x})=f(\sqrt{x}-1)$ is equivalent to the condition that $f(\sqrt{x})$ is a polynomial in $x+\sqrt{x}$. Indeed any nonzero polynomial $f(\sqrt{x})$ satisfying the first condition must have leading term of the form $a(\sqrt{x})^{2k}$. Subtracting $a(x+\sqrt{x})^k$ then gives another polynomial satisfying the same condition, but of strictly smaller degree. The claim therefore follows by induction on the degree of $f$. Next note that the image of $(x+\sqrt{x})^k$ under the linear map has leading term $x^k$. The result follows. $\Box$ 

\section{Algebras $D(q)$}
\begin{defi} Suppose $n \geq 4$ is an integer, $q(x)\in \mathbb{C} [x]$ is a polynomial of degree $n$. Define $\rho:=2q(-\tfrac{1}{2})$. Also define
\[s(x):=\frac{q(x)q(-x-1)}{x(x+1)(1+2x)^2}, \hspace{10mm}p(x):=\frac{-4q(x)q(-x-1)+\rho^2}{(1+2x)^2}.\]
Also, define $D(q)$ to be the subalgebra of $\tilde{T}(s)$ generated by the three elements
\[u := h^2,\hspace{10mm} v := a+b+\frac{2\rho}{1-4h^2},\hspace{10mm} w := (a-b)h-\frac{\rho}{1-4h^2}.\]
\end{defi}
\noindent
\begin{lemma}
\label{manyparts}
i) $\displaystyle{\hspace{5mm} s(x)=\frac{1}{4x(x + 1)}\bigg( -p(x) + \frac{\rho^2}{(1+2x)^2}\bigg),}$
\newline
ii) $p(-x)=p(x-1)$,
\newline
iii) $p(x)$ is a polynomial,
\newline
iv) $\{ p(\sqrt{x})\}$ and $\{ \sqrt{x} p(\sqrt{x})\}$ have degrees $n-2$ and $n-1$ respectively.
\end{lemma}

\noindent
Proof. The first two parts are easy to check. For iii), $p(x)$ is a polynomial because both the numerator and its first derivative give zero on substituting $x=-\tfrac{1}{2}$. For iv), assume that $q(x)$ has leading term $\alpha x^n$.  Clearly $p(\sqrt{x})$ has leading term $(-1)^{n+1} \alpha^2 x^{n-1}$, so $\{ \sqrt{x} p(\sqrt{x})\}$ has leading term $(-1)^{n+1} \alpha^{2} x^{n-1}$ also. In fact from ii) and the proof of Lemma \ref{iso} we know that $p(\sqrt{x})$ is a polynomial in $x + \sqrt{x}$, so must actually begin $(-1)^{n+1}\alpha^2 x^{n-1} + (-1)^{n+1}(n-1)\alpha^2 \sqrt{x} x^{n-2}$. It follows that $\{ p(\sqrt{x})\}$ has leading term $(-1)^{n+1}(n-1)\alpha^2 x^{n-2}$. $\Box$

\begin{theo}
\label{dqdef}
In $D(q)$ the following relations hold.
\[\left[u,v\right] = 2w + v,\hspace{5mm}\left[u,w\right]  =  2vu + w +\rho, \hspace{5mm}\left[v,w\right] =  -v^2 - \big\{ p(\sqrt{u})  \big\} ,\] 
\[w^2 - v^2 u - vw - \rho v - \big\{ \sqrt{u}p(\sqrt{u}) \big\} =0. \]
Moreover, $D(q)$ coincides with the algebra defined by this presentation.
\end{theo}

\noindent
Proof. We shall prove these relations in turn. Firstly,
\[
\left[u,v\right]  =  h^2 a + h^2 b - ah^2 - bh^2 = a(h+1)^2 + b(h-1)^2 - ah^2 - bh^2 =  2(a-b)h + a + b =  2w + v,\]
as required. Next,
\[
\left[u,w\right]  =  h^2ah - h^2 bh - ah^3 + bh^3 = a(h+1)^2 h - b(h-1)^2 h - ah^3 + bh^3 = 2(a+b)h^2 + (a-b)h \]\[=  2\left(v - \frac{2\rho}{1-4h^2}\right)h^2 + w + \frac{\rho}{1-4h^2} =  2vu + w + \rho.\] 
The final two relations are harder to see. To simplify the calculations we shall abbreviate $\displaystyle \frac{\rho}{1-4x^2}$ to $f(x)$. $\left[v,w\right] + v^2$ is equal to 
\[(a+b+2f(h)) ( (a-b)h - f(h)) -( (a-b)h - f(h))( a+b+2f(h))+(a+b+2f(h)) (a+b+2f(h)).\]
It is tedious but straightforward to bring all rational functions in $h$ to the right. This gives
\begin{eqnarray*}
\left[v,w\right] + v^2 &  = &  2(ba - ab)h + 2(ab + ba) + 4f(h)^2 \\
&& + a( 3f(h+1) + f(h) +2f(h+1)h - 2hf(h)) \\
&& + b( 3f(h-1) + f(h) -2f(h-1)h + 2hf(h)).\\
\end{eqnarray*}
Now, the coefficient of $a$ in this sum is zero. This follows from the identity
\[\frac{3}{1-4(x+1)^2 } + \frac{1}{1- 4x^2} + \frac{2x}{1-4(x+1)^2 } -\frac{2x}{1-4x^2 } = 0.\]
Similarly the coefficient of $b$ is zero. We also have
\vspace{5mm}
\newline
$2(ba - ab)h + 2(ab + ba) + 4f(h)^2$
\begin{eqnarray*}
& = & 2(h + 1)s(h) - 2(h - 1)s(h-1) + 4f(h)^2 \\
& = & \frac{1}{2h}\bigg( -p(h) + \frac{\rho^2}{(1+2h)^2}\bigg) - \frac{1}{2h}\bigg( -p(h-1) + \frac{\rho^2}{(1-2h)^2}\bigg) +\frac{4\rho^2}{(1-4h^2)^2} \\
& = & - \big\{ p(\sqrt{u})\big\} + \frac{\rho^2}{2h(1-4h^2)^2}\left( (1-2h)^2 - (1+2h)^2 + 4(2h) \right) \\
& = & - \big\{ p(\sqrt{u})\big\},
\end{eqnarray*}
as required. Here we used parts i) and ii) of Lemma \ref{manyparts} as well as Lemma \ref{root2}.
\par
Only the fourth relation remains to be proved. We have, by another tedious but easy calculation, that
\vspace{5mm}
\newline
$w^2 - v^2 u - vw - \rho v $ 
\begin{eqnarray*}
& = & ( (a-b)h-f(h) )^2 - ( a+b +2f(h) )^2 h^2 - ( a+b + 2f(h) )( (a-b)h - f(h) ) \\&&- ( a+ b+2f(h)) ( 1-4h^2 )f(h) \\
& = & -2ab(h^2 -h) -2ba(h^2 + h) + (1+4h^2)f(h)^2 \\&&-a(3f(h+1) + f(h) + 2f(h+1)h - 2f(h)h) h \\&&+b(3f(h-1) + f(h) - 2f(h-1)h + 2f(h)h) h \\
& = & -\frac{1}{2}\bigg( -p(h-1) + \frac{\rho^2}{(1-2h)^2}\bigg) -\frac{1}{2}\bigg( -p(h) + \frac{\rho^2}{(1+2h)^2}\bigg) +\frac{\rho^2(1+4h^2)}{(1-4h^2)^2} \\
& = & \big\{ \sqrt{u} p(\sqrt{u})\big\} -\frac{\rho^2}{2(1-4h^2)^2}( (1+2h)^2 + (1-2h)^2 - 2(1+4h^2)) \\
& = &\big\{ \sqrt{u} p(\sqrt{u})\big\}.
\end{eqnarray*}
This establishes that the fourth relation holds. 
\par To prove the final part of the theorem, we need to consider the algebra $A$ with generators $\tilde{u}$, $\tilde{v}$, $\tilde{w}$ defined by these four relations. It is clear from the relations that $A$ is spanned by the monomials of the form $\tilde{v}^i \tilde{w}^j \tilde{u}^k$ where $i,j,k\in \mathbb{N}$ and $j\leq 1$; indeed the four relations can be used respectively to deal with the `bad' submonomials $\tilde{u}\tilde{v}$, $\tilde{u}\tilde{w}$, $\tilde{w}\tilde{v}$, $\tilde{w}^2$. Thus the result will follow if we can prove that in $\tilde{T}(s)$ the elements $v^i w^j u^k$ where $i,j,k\in \mathbb{N}$ and $j\leq 1$ are linearly independent. We shall in fact prove the slightly stronger statement that, considering $\tilde{T}(s)$ as a right $\mathbb{C} (h)$-module, the elements $1,v,v^2 ,\cdots, w, vw, v^2w ,\cdots$ are independent over $\mathbb{C} (h)$. Indeed, suppose $k\geq 0$ and we have a relation
\[X:= c_0 (h) + \sum_{i=1}^{k}(v^i c_i (h) + v^{i-1} w d_i (h)) = 0,\]
where the $c_i (h)$ and the $d_i (h)$ are rational functions. If $k=0$ then obviously $c_0 (h)=0$, so assume $k\geq 1$. Now we have from the definition of $\tilde{T}(s)$ as a skew group algebra that $\tilde{T}(s)=\bigoplus_{i\in \mathbb{Z}} a^i \mathbb{C} (h)$. Moreover $(a^i \mathbb{C} (h) )(a^j \mathbb{C} (h))\subseteq a^{i+j} \mathbb{C} (h)$ for all $i,j\in \mathbb{Z}$, and $b=a^{-1}s(h-1)\in a^{-1} \mathbb{C} (h)$. It follows from the definition of $v$ and $w$ that the `$a^k$-term' of $X$ is $a^k (c_k (h) + hd_k(h))$. Hence $c_k (h) + hd_k (h) = 0$. Similarly from the $a^{-k}$-term we get $c_k (h) - hd_k (h) = 0$. Hence $c_k (h) = d_k (h) = 0$. The result therefore follows by induction on $k$. $\Box$   

\par
From now on we shall use these two descriptions of $D(q)$ interchangeably. We shall regard $D(q)$ as being a filtered algebra by putting $u$, $v$, $w$ in degrees $4$, $2n-4$, $2n-2$ respectively. 

\begin{rema}
\label{scalar}
If $\theta$ is a nonzero complex number, then $D(\theta q)\cong D(q)$ since the elements $u$, $\theta v$ and $\theta w$ of $D(q)$ can easily be seen to satisfy the relations of $D(\theta q)$.
\end{rema} 

The following is a simple consequence of the proof of Theorem \ref{dqdef}.

\begin{theo}
$D(q)$ is a deformation of the coordinate algebra of a Kleinian singularity of type $D_n$. $\Box$
\end{theo}

\begin{defi}
\label{const}
It follows from the proof of Theorem \ref{dqdef} that any $X\in D(q)$ can be written in the form $c_0 (u) + \sum_{i\geq 1}(v^i c_i (u) + v^{i-1} w d_i (u))$, where $c_i(x),d_i(x)\in \mathbb{C} [x]$ are uniquely determined polynomials. Define $c(X):=c_0(u)$ to be the constant term.
\end{defi}  

\noindent
One slightly subtle consequence of the proof of Theorem \ref{dqdef} is that if an element $X$ of $D(q)$ can be written in the form $c_0 (h) + \sum_{i=1}^{k}(v^i c_i (h) + v^{i-1} w d_i (h))$, where the $c_i (h)$ and the $d_i (h)$ are rational functions in $h$, then the $c_i (h)$ and the $d_i (h)$ must in fact be polynomials in $u$, so in particular $c(X) = c_0 (h)$. We shall use this idea in the proof of the following lemma. We shall require the following result in the proof Theorem \ref{theo2}.
\begin{lemma}
\label{commutation}
Suppose $t(x)\in \mathbb{C} [x]$. Then we have 
\[c\left( \left\{t\left(\sqrt{u}\right)\right\}v\right)=\frac{\rho}{1-4u}\left\{\frac{\left( 1-2\sqrt{u}\right) (  t\left( \sqrt{u} +1\right) - t\left( -\sqrt{u} -1\right) )}{2\left( \sqrt{u} + 1\right)}+2t\left( \sqrt{u} \right)\right\},\]
\[c\left( \left\{t\left(\sqrt{u}\right)\right\}w\right)=\frac{\rho}{1-4u}\left\{\frac{\sqrt{u}\left( 1-2\sqrt{u}\right) (  t\left( \sqrt{u} +1\right) - t\left( -\sqrt{u} -1\right) )}{2\left( \sqrt{u} + 1\right)}-t\left( \sqrt{u} \right)\right\}.\]
\end{lemma}

\noindent
(The contents of the braces are easily seen to be polynomials in $\sqrt{u}$. It is less obvious that the right hand sides of the above equalities are polynomials in $u$, but it is true and it can either be seen directly or deduced as a consequence of the following proof.)

\noindent
Proof. It follows immediately from the definitions of the elements $v$ and $w$ that
\[ah+bh=vh-\frac{2\rho h}{1-4h^2},\hspace{10mm} ah-bh=w+\frac{\rho}{1-4h^2}.\]
It follows easily that
\[a=\frac{1}{2}v + \frac{1}{2}w\frac{1}{h}+\frac{\rho}{2h(1+2h)},\hspace{10mm} b=\frac{1}{2}v - \frac{1}{2}w\frac{1}{h}-\frac{\rho}{2h(1-2h)}.\]
Now we get
\begin{eqnarray*}
c\left( \left\{t\left(\sqrt{u}\right)\right\}v\right) & = & c\left( \frac{t(h)-t(-h)}{2h}\left( a + b +\frac{2\rho}{1-4h^2}\right)\right) \\
 &= & c\left( a\frac{t(h+1)-t(-h-1)}{2(h+1)}+b\frac{t(h-1)-t(-h+1)}{2(h-1)}+\frac{2\rho \left\{t( \sqrt{u}) \right\}}{1-4h^2}\right) \\
 &= & \frac{\rho}{2h(1+2h)}\frac{t(h+1)-t(-h-1)}{2(h+1)}-\frac{\rho}{2h(1-2h)}\frac{t(h-1)-t(-h+1)}{2(h-1)} \\
 &&+\frac{2\rho \left\{t( \sqrt{u}) \right\}}{1-4h^2} \\
& = & \frac{\rho}{1-4h^2}\left( \frac{\frac{(1-2h)(t(h+1)-t(-h-1))}{2(h+1)}-\frac{(1+2h)(t(-h+1)-t(h-1))}{2(-h+1)}}{2h}+2\left\{ t( \sqrt{u})\right\} \right) \\
& = & \frac{\rho}{1-4u}\left\{\frac{\left( 1-2\sqrt{u}\right) (  t\left( \sqrt{u} +1\right) - t\left( -\sqrt{u} -1\right) )}{2\left( \sqrt{u} + 1\right)}+2t\left( \sqrt{u} \right)\right\},
\end{eqnarray*}
\noindent
as required. The proof of the second equality is similar. $\Box$

\section{Levy's Presentations}

In this section we shall relate the algebras $D(q)$ to the algebras $D(Q,\rho)$ constructed in \cite{levy} and derive some consequences. 
\begin{defi}
Suppose that $n\geq 4$, $Q(X)$ is a polynomial of degree $n-1$ and $\rho\in\mathbb{C}$. Let $P(X)$ be the unique polynomial satisfying the condition 
\[Q(-X(X+1))+(X+1)P(-X(X+1)) \text{ is a polynomial in } X^2.\]
Define $D(Q, \rho)$ to be the algebra with generators $x$, $y$, $z$ satisfying the relations
\[\left[ x,y\right] = 2z,\hspace{5mm} \left[x,z\right]=-2xy+2z+\rho,\hspace{5mm}\left[y,z\right]=y^2+P(x),\hspace{5mm}Q(x)+xy^2+z^2-2zy-\rho y=0.\]
$D(Q,\rho)$ is made into a filtered algebra by putting $x$, $y$, $z$ in degrees $4$, $2n-4$, $2n-2$ respectively. 
\end{defi}

It is shown in \cite{levy} that $D(Q,\rho)$ is a deformation of a Kleinian singularity of type $D_n$, that \emph{every} deformation of a Kleinian singularity of type $D_n$ that is not commutative is isomorphic to some $D(Q, \rho)$, and, for $n\geq 5$, that $D(Q,\rho)\cong D(Q',\rho ')$ if and only if there is a nonzero scalar $\theta$ such that $Q' = \theta^2 Q$ and $\rho' = \theta \rho$ (This last part is part a) of Theorem 2.22). The `if' part of this last statement is obvious; indeed the elements $x$, $\theta y$, $\theta z$ of $D(Q, \rho)$ clearly satisfy the relations for $D(\theta^2 Q,  \theta\rho)$. 
\par
We shall now explain the connection between the families $D(q)$ and $D(Q, \rho)$. Define the elements $u$, $v$ of $w$ of $D(Q, \rho)$ by 
\[u:=-x+\tfrac{1}{4},\hspace{10mm} v:=y, \hspace{10mm} w:=-z-\tfrac{1}{2}y.\]
It is easy to check that $u$, $v$, $w$ satisfy the relations
\[\left[u,v\right]=2w+v, \hspace{5mm} \left[u,w\right]=2vu+w+\rho,\hspace{5mm}\left[ v,w\right]=-v^2-P(-u+\tfrac{1}{4}),\]
\[w^2-v^2u-vw-\rho v+Q(-u+\tfrac{1}{4})+\tfrac{1}{2}P(-u+\tfrac{1}{4})=0.\]
Now define a polynomial $p$ in a formal square root of an indeterminate $X$ by \[p(\sqrt{X}):=-Q(-X+\tfrac{1}{4})+(\sqrt{X}-\tfrac{1}{2})P(-X+\tfrac{1}{4}).\]
It is clear that
\[\left\{p(\sqrt{u})\right\}=P(-u+\tfrac{1}{4}),\hspace{10mm}\left\{ \sqrt{u}p(\sqrt{u})\right\}=-Q(-u+\tfrac{1}{4})-\tfrac{1}{2}P(-u+\tfrac{1}{4}).\]
Now substituting $-\sqrt{u}-\tfrac{1}{2}$ for $X$ in the condition connecting $P$ and $Q$ we deduce that $p(\sqrt{u})$ is a polynomial in $(\sqrt{u}+\tfrac{1}{2})^2$, hence is a polynomial in $u+\sqrt{u}$. The same is therefore true of the polynomial $\rho^2 - (1+2\sqrt{u})^2 p(\sqrt{u})$. It follows that we can write this last polynomial in the form $4q(\sqrt{u})q(-\sqrt{u}-1)$, and it is clear that $\rho=\pm 2q(-\tfrac{1}{2})$. If $\rho = -2q(-\tfrac{1}{2})$ we simply replace $q$ by its negative. In either case we have shown that $D(Q, \rho)$ is isomorphic to $D(q)$ for some $q$.
\par
We can also carry out the process described in reverse. This time let us start with a polynomial $q$ of degree $n$. Let $u$, $v$, $w$, $\rho$, $p$ be the usual notation associated with the algebra $D(q)$. We may define elements $x$, $y$, $z$ of $D(q)$ by using the same equations we used to define $u$, $v$, $w$ before. We may also define polynomials $P$ and $Q$ by inverting the process used to define $p$ before; in other words we define $P$ and $Q$ by
\[Q(-u+\tfrac{1}{4})=\left\{(-\sqrt{u}-\tfrac{1}{2})p(\sqrt{u})\right\},\hspace{10mm}P(-u+\tfrac{1}{4})=\left\{p(\sqrt{u})\right\}.\]
As before it is a simple matter to check that $x$, $y$, $z$ satisfy the required relations and that the polynomials $P$ and $Q$ satisfy the desired condition.   
\begin{theo}
\label{theo1}
If $n\geq 4$ every deformation of a Kleinian singularity of type $D_n$ that is not commutative is isomorphic as a filtered algebra to $D(q)$ for some $q$ of degree $n$. If $n\geq 5$ and $q$, $q'$ are two polynomials of degree $n$, then $D(q)\cong D(q')$ as filtered algebras if and only if there exists a nonzero scalar $\theta$ such that $q'(x)q'(-x-1)=\theta^2 q(x)q(-x-1)$.
\end{theo}
\noindent
Proof. The first part follows immediately from the above discussion. The `if' part of the second statement follows from Remark  \ref{scalar} and the definition of $D(q)$. For the `only if' part, suppose that $D(q)\cong D(q')$. Then, the corresponding result for the family $D(Q, \rho)$ implies, by the above discussion, that there exists some nonzero $\theta$ such that
\[\left\{(-\sqrt{x}-\tfrac{1}{2})p'(\sqrt{x})\right\}=\theta^2 \left\{(-\sqrt{x}-\tfrac{1}{2})p(\sqrt{x})\right\},\hspace{5mm} \rho' = \theta \rho,\]
where $p,\rho$ are as usual for $D(q)$ and $p',\rho '$ are the corresponding things for $D(q')$. Next we apply Lemma \ref{iso} to deduce that $p'(x)$=$\theta^2 p(x)$. Thus $4q'(x)q'(-x-1)-\rho'^2 = \theta^2(4q(x)q(-x-1)-\rho^2)$. The result follows. $\Box$

We will see in the proof of Theorem \ref{theo3} the necessary and sufficient condition for two $D(q)$ to be isomorphic in the case $n=4$.

\begin{prop}
\label{mainprop}
Suppose that $A$ is a filtered algebra with filtration denoted $F_0\subseteq F_1 \cdots $ whose associated graded algebra is isomorphic to the coordinate algebra of a Kleinian singularity of type $D_{n}$. Suppose that $x\in F_4$ and $y\in F_{2n-4}$ are non-commuting elements. Define elements $z$ and $\sigma$ by the equations $\left[ x,y\right] =2z+y$ and $\left[ x,z\right] = 2yx + z + \sigma$. Suppose further that $\sigma $ is a scalar multiple of $1$. Then for a suitable choice of $q'$ there is an isomorphism $A\cong D(q')$ identifying $x$, $y$ and $z$ with the standard generators $u'$, $v'$ and $w'$ of $D(q')$.
\end{prop}

\noindent
Proof. We may assume that $A=D(q)$ for some polynomial $q$ of degree $n$. Let $u$, $v$, $w$, $p$, $\rho$ be the standard notation associated with $D(q)$. 
\par
Case 1: $n\geq 5$. Since $x$ and $y$ are elements of $F_4$ and $F_{2n-4}$ respectively we can write $x=\lambda u + a$, $y=\mu v + r(u)$ where $a, \lambda$, $\mu\in \mathbb{C}$ and $r$ is a polynomial. Note that we cannot do this for the case $n=4$ because then $x$ could have a coefficient of $v$. Since $\left[ x,y\right]\neq 0$ we must have $\lambda \mu \neq 0$, so by replacing $v$ by a nonzero scalar multiple we can assume without loss that $\mu  =1$. Now $2z + y =\left[ x,y\right] = \lambda \left[ u,v\right]=\lambda (2w+v)$. Therefore $z=\lambda w + \tfrac{\lambda}{2} v - \tfrac{1}{2} y = \lambda w + \tfrac{\lambda -1}{2}v - \tfrac{1}{2}r(u)$. Next, from $\left[ x,z\right] = 2yx+ z + \sigma$ we get
\begin{eqnarray*}
\sigma & = & \left[ \lambda u + a,\lambda w + \tfrac{\lambda - 1}{2} v - \tfrac{1}{2}r(u)\right] - 2(v+r(u))(\lambda u  + a) - \lambda w - \tfrac{\lambda - 1}{2}v + \tfrac{1}{2}r(u) \\
& = & \lambda^2 \left[ u,w\right] + \tfrac{\lambda (\lambda -1)}{2}\left[ u,v\right] - 2\lambda vu -2av - 2\lambda r(u)u - 2ar(u) - \lambda w - \tfrac{\lambda - 1}{2}v + \tfrac{1}{2}r(u) \\
& = & \lambda^2 (2vu + w + \rho) + \tfrac{\lambda (\lambda -1)}{2} (2w + v)\\
&& - 2\lambda vu -2av - 2\lambda r(u)u - 2ar(u) - \lambda w - \tfrac{\lambda - 1}{2}v + \tfrac{1}{2}r(u)
\end{eqnarray*}
Since $D(q)$ has basis the monomials $v^i w^j u^k$ where $j\leq 1$ we can compare coefficients. The coefficient of $vu$ gives $2\lambda^2 = 2\lambda$, so $\lambda = 1$. The coefficient of $v$ then tells us that $a=0$. Finally, comparing terms not containing $v$ or $w$ we get $\sigma = \rho + r(u)(-2u + \tfrac{1}{2})$. It follows that $r(u)=0$, for otherwise the right hand side would not be constant. We have therefore proved that $(x,y)=(u,v)$. The result follows. 
\par
Case 2: $n=4$. Since $x$, $y$ are elements of $F_4$ we can write $x=au+bv+e$ and $y=cu+dv+f$, where $a,b,c,d,e,f\in\mathbb{C}$. Now $\left[ x,y\right] = (ad-bc)\left[ u,v\right]$. Since $x$ and $y$ don't commute we must have $ad-bc \neq 0$. In fact by replacing $v$ by a suitable nonzero scalar multiple we can assume without loss that $ad - bc = 1$. We shall now show that these reductions only leave three possibilities for the pair $(x,y)$. We have $2z + y = \left[ x,y\right] = \left[ u,v\right] = 2w + v$. Therefore $z=w+\frac{1}{2}v-\frac{1}{2}y = w+\frac{1-d}{2}v-\frac{c}{2}u-\frac{f}{2}$. Next, from $\sigma = \left[ x,z\right] - 2yx - z$ we get
\begin{eqnarray*}
\sigma & = & \left[ au + bv + e,w+\tfrac{1-d}{2}v-\tfrac{c}{2}u-\tfrac{f}{2}\right] - 2(cu+dv+f)(au+bv+e)\\
&& -w-\tfrac{1-d}{2}v+\tfrac{c}{2}+\tfrac{f}{2} \\
& = & a\left[ u,w\right] + \tfrac{a(1-d)+bc}{2}\left[ u,v\right]+b\left[ v,w\right]- 2acu^2-2bcuv-2ceu-2advu-2bdv^2-2dev \\
&&-2afu-2bfv-2ef-w-\tfrac{1-d}{2}v+\tfrac{c}{2}+\tfrac{f}{2} \\
& = & a(2vu+w+\rho)+\tfrac{a-1}{2}(2w+v)+b(-v^2 - \left\{p(\sqrt{u})\right\})-2acu^2-2bc(vu+2w+v)\\
&& -2ceu-2advu - 2bdv^2 - 2dev - 2afu - 2bfv - 2ef - w-\tfrac{1-d}{2}v+\tfrac{c}{2}u+\tfrac{f}{2}.
\end{eqnarray*}  
Since $D(q)$ has basis the monomials $v^i w^j u^k$ where $j\leq 1$ we can compare coefficients. From the coefficient of $v^2$ we deduce $b(1+2d)=0$, so either $b=0$ or $d=-\tfrac{1}{2}$. From the coefficient of $vu$ we get $2a-2bc-2ad=0$, hence $a-2ad+1=0$ from $ad-bc=1$. (This is in fact the equation one gets by comparing the coefficients of $w$ too). The coefficient of $v$ gives 
\[0=\tfrac{a}{2}-2ad+1-2de-2bf+\tfrac{d}{2}, \hspace{10mm} (*)\] as can readily be verified. Finally comparing terms not involving $v$ or $w$ we get
\[\sigma = a\rho - b\left\{p(\sqrt{u})\right\}-2acu^2-2ceu-2afu-2ef+\tfrac{c}{2}u+\tfrac{f}{2}.\hspace{10mm} (**)\]
Let us suppose first that $b=0$. Then from the coefficient of $u^2$ in $(**)$ we get $ac=0$. Since $b=0$ and $ad-bc=1$ we cannot have $a=0$ so we must have $c=0$. The coefficient of $u$ in $(**)$ gives $f=0$. Now, since $ad-bc=1$ and $b=0$ we deduce $ad=1$. Since $a-2ad+1=0$ we get $a=d=1$. The only letter from $a$-$f$ we have not yet determined is $e$, but we can now deduce that $e=0$ from $(*)$. This shows that $(x,y)=(u,v)$. Next, let us suppose that $b\neq 0$. Then $d=-\tfrac{1}{2}$. From $a-2ad+1=0$ we get $a=-\tfrac{1}{2}$. With these values, equation $(*)$ now reads $e=2bf$. Since $ad-bc=1$ and $a=d=-\tfrac{1}{2}$ we get $c=-\tfrac{3}{4b}$. Substituting the values for $a$, $d$ and these expressions for $e$ and $c$ into $(**)$ gives
\[\sigma = -\tfrac{1}{2}\rho - b\left\{ p( \sqrt{u})\right\} - \tfrac{3}{4b}u^2 +fu - 4bf^2 - \tfrac{3}{8b}u+\tfrac{f}{2}.\]
Now suppose $\theta$ is the coefficient of $u^2$ in $\left\{ p( \sqrt{u})\right\}$. From Lemma \ref{manyparts} part iv) we know that $\theta \neq 0$. We get $b\theta +\tfrac{3}{4b}=0$, so that $b$ is one of the two distinct roots of $X^2 + \tfrac{3}{4\theta}=0$. Whichever choice for $b$ we make, the value of $f$, and hence $e$, is then determined by considering the coefficient of $u$ in the above. We have therefore shown that there are at most three possibilities for the pair $(x, y)$. But it follows from \cite{levy} (Theorem 3.6(a)) that there are precisely three non-trivially different isomorphims $D(Q_i ,\rho_i)\rightarrow D(q)$ for $i=1,2,3$. By non-trivially different we mean that no two of the pairs $(x_i ,y_i)$ are such that $x_i=x_j$ and $y_i$ is a scalar multiple of $y_j$ (here $x_i,y_i,z_i$ are the images in $D(q)$ of the standard generators for $D(Q_i ,\rho_i)$). It follows from the discussion at the beginning of this section that there are exactly three possibilities for the pair $(x,y)$ and that these must have the desired property. $\Box$

\section{Algebras $\Pi^{\boldsymbol{\lambda}}$ and $\mathcal{O}^{\boldsymbol{\lambda}}$}
We begin by fixing some notation. Define $Q$ to be the quiver

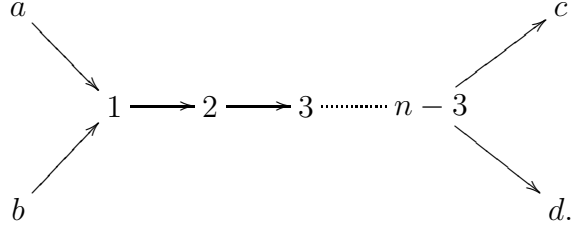
\[
\xymatrix{ a\ar[rd] & & & & & c \\ & 1\ar[r] & 2\ar[r] & 3\ar@{.}[r] & n-3\ar[ur]\ar[dr] & \\ b\ar[ur] & & & & & d.}
\]
\noindent
It is conventional to denote the affine vertex by $0$ but we shall use $a$ in order to emphasize the symmetry in $a$, $b$, $c$ and $d$. 
Let $I$ denote the set of vertices of $Q$, and $\mathbb{C}^I$ the vector space with basis the set of vertices. We shall use bold script to denote elements of $\mathbb{C}^I$. Given two elements $\boldsymbol{\alpha}$, $\boldsymbol{\beta}$ of $\mathbb{C}^I$, define their dot product by $\boldsymbol{\alpha}\cdot\boldsymbol{\beta}=\sum_I \boldsymbol{\alpha}^i \boldsymbol{\beta}^i$, where  $\boldsymbol{\alpha}^i$ denotes the coefficient of $i$ in $\boldsymbol{\alpha}$. Define $\boldsymbol{\delta}$ by

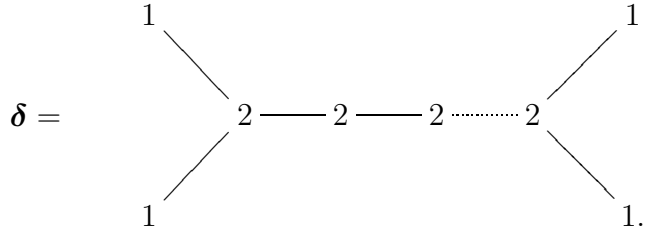
\[
\xymatrix{ & 1\ar@{-}[rd] & & & & & 1 \\ \boldsymbol{\delta} = && 2\ar@{-}[r] & 2\ar@{-}[r] & 2\ar@{.}[r] & 2\ar@{-}[ur]\ar@{-}[dr] & \\ & 1\ar@{-}[ur] & & & & & 1.}
\]
Also fix a vector 
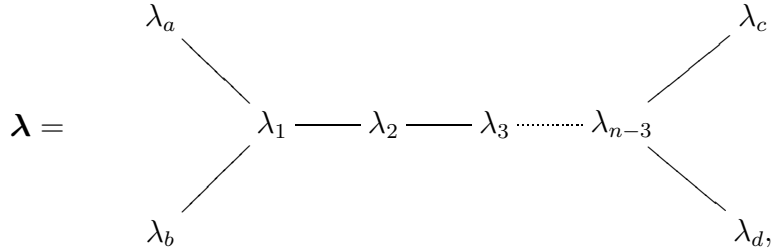
\[
\xymatrix{ & \lambda_a\ar@{-}[rd] & & & & & \lambda_c \\ \boldsymbol{\lambda} = && \lambda_1 \ar@{-}[r] & \lambda_2 \ar@{-}[r] & \lambda_3 \ar@{.}[r] & \lambda_{n-3} \ar@{-}[ur]\ar@{-}[dr] & \\ & \lambda_b \ar@{-}[ur] & & & & & \lambda_d,}
\]
which we shall assume satisfies $\boldsymbol{\lambda}\cdot\boldsymbol{\delta} = 1$. 

\begin{defi}
Define complex numbers $\mu_0$,...,$\mu_{n-1}$ by $\mu_0 := \tfrac{1}{2}\lambda_a - \tfrac{1}{2}\lambda_b$, $\mu_1 := \tfrac{1}{2}\lambda_a + \tfrac{1}{2}\lambda_b$, $\mu_2 := \tfrac{1}{2}\lambda_a + \tfrac{1}{2}\lambda_b + \lambda_1$, $\cdots$, $\mu_{n-2} := \tfrac{1}{2}\lambda_a + \tfrac{1}{2}\lambda_b + \lambda_1 + \cdots +\lambda_{n-3}$, $\mu_{n-1} := \tfrac{1}{2}\lambda_a + \tfrac{1}{2}\lambda_b + \lambda_1 + \cdots +\lambda_{n-3} + \lambda_c$. Define also $\gamma := \tfrac{1}{2}\lambda_c - \tfrac{1}{2}\lambda_d$.
\end{defi}

\begin{defi}
If $i$, $j\in I$, let $\lambda_{ij}$ or sometimes $\lambda_{i,j}$ denote the complex number $\sum \lambda_k$ where the sum is over all elements $k$ of $I$ lying on the shortest path connecting $i$ and $j$ on the underlying graph of $Q$ (this includes both $i$ and $j$). 
\end{defi}
\begin{defi}
Let $\overline{Q}$ denote the double quiver of $Q$, that is the quiver obtained by adjoining for each arrow $a: i\rightarrow j$ a reverse arrow $a^{\ast}: j \rightarrow i$. Let $\mathbb{C} \overline{Q}$ denote the path algebra of $\overline{Q}$. Given a non-empty word $i_1 i_2 ... i_k$ in the elements of $I$, define $\langle i_1 i_2 ... i_k \rangle\in \mathbb{C} \overline{Q}$ to be the unique shortest path starting at $i_1$, ending at $i_k$ and visiting the vertices in the correct order. Thus, for example $\langle aba \rangle$ is a path of length $4$, $\langle acba \rangle$ is a path of length $2n - 2$, while $\langle d \rangle$ denotes the vertex idempotent at the vertex $d$. We shall extend this by linearity to define an element $\langle X \rangle\in \mathbb{C} \overline{Q}$ for any element $X$ of the free associative $\mathbb{C}$-algebra (without a $1$) generated by the set $I$. Thus, for example, $\langle a(a+b)a \rangle$ is the vertex idempotent at $a$ added to the path $\langle aba \rangle$ of length $4$. 
\end{defi}

There is ambiguity in this notation because we are using numbers to denote vertices. For example $\langle 2a \rangle$ could mean twice the vertex idempotent at $a$, or the path from the vertex $2$ to the vertex $a$. Therefore we shall always underline all the vertices with numbers for names within $\langle$ and $\rangle$.

The following definitions first appeared in \cite{CBH}.

\begin{defi}
The deformed preprojective algebra $\Pi^{\boldsymbol{\lambda}}$ is defined to be the quotient of $\mathbb{C} \overline{Q}$ given by the relation $\sum \left[ \alpha, \alpha^{\ast} \right] = \sum \lambda_i \langle i \rangle$, where the first sum is over all arrows in $Q$ and the second is over the set $I$. We define $\mathcal{O}^{\boldsymbol{\lambda}}$ to be the algebra $\langle a \rangle \Pi^{\boldsymbol{\lambda}} \langle a \rangle$.
\end{defi}
\begin{rema}
\label{orientation}
The construction of $\Pi^{\boldsymbol{\lambda}}$ is independent of the orientation of the graph $Q$. Indeed if we replace an arrow $\alpha : i\rightarrow j$ in $Q$ by a reverse arrow $\beta :j\rightarrow i$ then the new algebra is isomorphic to the old one by the map sending $\alpha$ to $\beta^{\ast}$ and $\alpha^{\ast}$ to $-\beta$. 
\end{rema}
\begin{rema}
The identity element of $\mathcal{O}^{\boldsymbol{\lambda}}$ is $\langle a\rangle$, so when we are considering elements of $\mathcal{O}^{\boldsymbol{\lambda}}$ we shall write $\alpha$ rather than $\alpha \langle a\rangle$ for any complex number $\alpha$.
\end{rema}
\begin{rema}
If $n \geq 5$ the relation defining $\Pi^{\boldsymbol{\lambda}}$ implies (by premultiplying or postmultiplying by the vertex idempotents) the following relations:
\begin{eqnarray*}
\langle a\hspace{1mm} \underline{1} \hspace{1mm} a  \rangle & = & \lambda_a \langle a \rangle\\
\langle b \hspace{1mm}\underline{1} \hspace{1mm} b  \rangle & = & \lambda_b \langle b \rangle\\
\langle \underline{1}\hspace{1mm}\underline{2}\hspace{1mm}\underline{1} \rangle - \langle \underline{1} \hspace{1mm} a \hspace{1mm} \underline{1} \rangle - \langle \underline{1} \hspace{1mm} b \hspace{1mm} \underline{1} \rangle & = & \lambda_1 \langle \underline{1} \rangle \\
\langle \underline{k} \hspace{1mm} \underline{k+1} \hspace{1mm} \underline{k} \rangle - \langle \underline{k} \hspace{1mm} \underline{k-1} \hspace{1mm} \underline{k} \rangle & = & \lambda_k \langle \underline{k} \rangle \hspace{10mm}(1<k<n-3) \\
\langle \underline{n-3} \hspace{1mm} c \hspace{1mm} \underline{n-3}\rangle + \langle \underline{n-3} \hspace{1mm} d \hspace{1mm} \underline{n-3}\rangle - \langle \underline{n-3} \hspace{1mm} \underline{n-4} \hspace{1mm} \underline{n-3}\rangle & = & \lambda_{n-3} \langle \underline{n-3} \rangle \\
\langle c\hspace{1mm} \underline{n-3} \hspace{1mm} c  \rangle & = & -\lambda_c \langle c \rangle\\
\langle d \,\underline{n-3} \, d  \rangle & = & -\lambda_d \langle d \rangle.
\end{eqnarray*}
\noindent
Clearly, one can write down similar relations for the case $n=4$.
\end{rema}
\begin{defi}
The algebra $\mathcal{O}^{\boldsymbol{\lambda}}$ is naturally filtered by the length of paths. We shall denote the filtration of $\mathcal{O}^{\boldsymbol{\lambda}}$ by $F_0 \subseteq F_1 \subseteq F_2 \subseteq \cdots$. 
\end{defi}
\par
Crawley-Boevey and Holland do not insist that $\boldsymbol{\lambda}\cdot\boldsymbol{\delta}=1$ in the definition of $\mathcal{O}^{\boldsymbol{\lambda}}$. However it is clear that for any nonzero scalar $\theta$ we have $\mathcal{O}^{\boldsymbol{\lambda}}\cong\mathcal{O}^{\theta\boldsymbol{\lambda}}$. Also, since we are only interested in deformations that are not commutative, and since $\mathcal{O}^{\boldsymbol{\mu}}$ is commutative if and only if $\boldsymbol{\mu}\cdot\boldsymbol{\delta}=0$ (Theorem 0.4 (1) of \cite{CBH}), we do not lose anything by making this assumption. 
\par
The proof of the following can be found in \cite{CBH}.
\begin{theo}
$\mathcal{O}^{\boldsymbol{\lambda}}$ is a deformation of the coordinate algebra of a Kleinian singularity of type $D_{n}$.
\end{theo} 

\section{Proof of the Main Theorem}

\begin{theo}
\label{theo2}
There is an isomorphism $\mathcal{O}^{\boldsymbol{\lambda}}\cong D(q)$ of filtered algebras where
\[q(x) :=  \prod_{i=0}^{n-1}(x + \mu_i ).\]
\end{theo}

\noindent
The proof of the theorem is fairly long, and shall be broken down into several small lemmas. Notice that $\rho:=2q(-\tfrac{1}{2})=2\gamma\prod_{i=0}^{n-2}\left(\mu_i - \tfrac{1}{2} \right)$. We shall now define elements $u$, $v$, $w$ of $\mathcal{O}^{\boldsymbol{\lambda}}$ that shall turn out to correspond under the isomorphism to the elements of $D(q)$ which were also called $u$, $v$, $w$. This abuse of notation should cause no confusion. 

\begin{defi}
The elements $u$, $v$, $w$ of $\mathcal{O}^{\boldsymbol{\lambda}}$ are defined as follows.
\[u :=\langle aba\rangle + \mu_0^2, \hspace{3mm}v := \langle a(c-d)a \rangle + \big\{ h(\sqrt{u}) \big\},\hspace{3mm}w  :=  \langle a(c-d)ba \rangle +\mu_0 v + \big\{ (\sqrt{u}-\mu_0)h(\sqrt{u})\big\},\]
\noindent
where 
\[h(x):=\frac{2\gamma\displaystyle{\prod_{i=0}^{n-2}}(x+\mu_i )-\rho}{2x+1}.\]
\end{defi}
\noindent
Notice that $h(x)$ is a polynomial in $x$ because by definition of $\rho$ the numerator of this fraction evaluated at $-\frac{1}{2}$ gives $0$.
We shall begin the proof of Theorem \ref{theo2} with a lemma.
\begin{lemma}
\label{kink}
Suppose $1\leq k\leq n-3$. Then in $\Pi^{\boldsymbol{\lambda}}$ we have 
\begin{eqnarray*}
\langle \underline{k} \, \underline{k+1} \, \underline{1} \rangle & = & \langle \underline{k} \,(a + b)\,\underline{1}\rangle + \lambda_{1k}\langle \underline{k}\,\underline{1}\rangle, \\
\langle \underline{k}\, \underline{k-1}\, \underline{n-3} \rangle & = & \langle \underline{k} \,(c+d)\,\underline{n-3}\rangle - \lambda_{k,n-3}\langle \underline{k}\,\underline{n-3}\rangle,
\end{eqnarray*}
where $\underline{0}$ and $\underline{n-2}$ are interpreted as $(a+b)$ and $(c+d)$ respectively.
\end{lemma}
\noindent
Proof. By symmetry it is only necessary to prove the first of these. If $k=1$ then the statement is simply the relation for $\langle \underline{1}\,\underline{2}\,\underline{1}\rangle$. If $k\geq 2$ then we use the relation for $\langle \underline{k}\,\underline{k+1}\,\underline{k}\rangle$ and induction to get
\[\langle\underline{k}\,\underline{k+1}\,\underline{1}\rangle = \langle \underline{k}\,\underline{k-1}\,\underline{k}\,\underline{1}\rangle + \lambda_k \langle\underline{k}\,\underline{1}\rangle = \langle \underline{k} \,(a + b)\,\underline{1}\rangle + \lambda_{1,k-1}\langle \underline{k}\,\underline{1}\rangle + \lambda_{k}\langle \underline{k}\,\underline{1}\rangle,\]
as required. $\Box$

\begin{lemma}
\label{zigzag}
In $\Pi^{\boldsymbol{\lambda}}$ we have the following relations.
\begin{eqnarray*}
\langle abc \rangle = \langle adc \rangle - \lambda_{ac} \langle ac \rangle, && \langle cda \rangle = \langle cba \rangle + \lambda_{ac} \langle ca \rangle, \\
\langle abd \rangle = \langle acd \rangle - \lambda_{ad} \langle ad \rangle, &&\langle cdb \rangle = \langle cab \rangle + \lambda_{bc} \langle cb \rangle, \\
\langle bac \rangle = \langle bdc \rangle - \lambda_{bc} \langle bc \rangle, &&\langle dca \rangle = \langle dba \rangle + \lambda_{ad} \langle da \rangle, \\
\langle bad \rangle = \langle bcd \rangle - \lambda_{bd} \langle bd \rangle, &&\langle dcb \rangle = \langle dab \rangle + \lambda_{bd} \langle db \rangle. 
\end{eqnarray*}
\end{lemma}
\noindent
Proof. We shall only prove the first of these, the rest being similar. We have
\[ \langle abc \rangle = \langle a \,\underline{1}\, (a+b)\, \underline{n-3}\, c\rangle - \langle a \,\underline{1} \,ac\rangle = \langle a (c+d) \,\underline{n-3}\, c\rangle - \lambda_{a,n-3}\langle ac\rangle,\]
by applying Lemma \ref{kink} to $\langle \underline{1} \,(a+b)\,\underline{n-3}\rangle$ and using the relation for $\langle a\,\underline{1}\,a\rangle$. Finally applying the relation for $\langle c\,\underline{n-3}\,c\rangle$ we get the result. $\Box$

\begin{lemma}
\label{recur1}
If $1 \leq k \leq n-2$ then in $\mathcal{O}^{\boldsymbol{\lambda}}$ we have
\[ \langle a \, \underline{k} \, a \rangle = \left\{ \prod_{i=0}^{k}(\sqrt{u} + \mu_i) \right\},\hspace{5mm} \langle a \, \underline{k} \, ba \rangle = \left\{ ( \sqrt{u} -\mu_0 )\prod_{i=0}^{k}(\sqrt{u} + \mu_i) \right\},\]
\noindent
where $\underline{n-2}$ is interpreted as $(c+d)$. \end{lemma}

\noindent
Proof. The proof is by induction on $k$. For $k=1$ we have
\[
\left\{ \prod_{i=0}^1 ( \sqrt{u} + \mu_i ) \right\}  =  \left\{ ( \sqrt{u} +\mu_0 ) ( \sqrt{u} + \mu_1 ) \right\} =  \mu_0 + \mu_1 =  \lambda_a =  \langle a \, \underline{1} \, a\rangle ,\]
\[
\left\{ ( \sqrt{u} - \mu_0 )\prod_{i=0}^1 ( \sqrt{u} + \mu_i ) \right\} =  \left\{ ( \sqrt{u} -\mu_0 )( \sqrt{u} +\mu_0 ) ( \sqrt{u} + \mu_1 ) \right\} = u-\mu_0^2  =  \langle a \, \underline{1} \, ba\rangle. \]
\noindent
Now we suppose that $k\geq 2$. Then we have
\begin{eqnarray*}
\langle a \, \underline{k} \, a\rangle & = & \langle a \, \underline{k-1} \, \underline{k} \,\underline{1} \, a \rangle \\
& =& \langle a \,\underline{k-1}(a+b)\underline{1}\, a\rangle + \lambda_{1,k-1}\langle a\,\underline{k-1}\, a\rangle\hspace{5mm}\text{applying Lemma \ref{kink} to }\langle \underline{k-1} \, \underline{k} \,\underline{1} \rangle \\
& = & \lambda_{a,k-1}\langle a \, \underline{k-1} \, a\rangle + \langle a \, \underline{k-1} \, ba\rangle \\
& = & \left\{ ( \lambda_{a,k-1} + \sqrt{u} - \mu_0 ) \prod_{i=0}^{k-1}( \sqrt{u} + \mu_{i}) \right\} \\
& = & \left\{ \prod_{i=0}^{k}( \sqrt{u} + \mu_i ) \right\}.
\end{eqnarray*}
\noindent
Similarly we have
\begin{eqnarray*}
\langle a \, \underline{k} \, ba \rangle & = & \lambda_{b, k-1}\langle a \, \underline{k-1} \,ba\rangle + \langle a \, \underline{k-1} \, aba \rangle \\
& = & \lambda_{b, k-1}\langle a \, \underline{k-1} \,ba\rangle + \langle a \, \underline{k-1} \, a \rangle \langle aba\rangle \\
& = & \lambda_{b, k-1}\langle a \, \underline{k-1} \,ba\rangle + \langle a \, \underline{k-1} \, a \rangle(u- \mu_0^2 ) \\
& = & \left\{ ( \sqrt{u} - \mu_0 )( \lambda_{b,k-1} + \sqrt{u} + \mu_0 )\prod_{i=0}^{k-1} (\sqrt{u}+\mu_i)\right\} \\
& = & \left\{ ( \sqrt{u} -\mu_0 )\prod_{i=0}^{k}(\sqrt{u} + \mu_i) \right\}.\hspace{5mm}\Box
\end{eqnarray*}

\par
Our method of proof of Theorem \ref{theo2} shall be to show that the elements $u$, $v$ satisfy the conditions (for $x$ and $y$) of Proposition \ref{mainprop}. We shall start this process with the following.  

\begin{lemma}
\label{relat1}
In $\mathcal{O}^{\boldsymbol{\lambda}}$ we have $0 \neq \left[u,v\right] = 2w + v$.
\end{lemma}

\noindent
Proof. We first claim that any path starting and ending at $a$ and not visiting either of the vertices $c$ or $d$ is a polynomial in $u$. Indeed suppose $k$ is the highest vertex number visited by the path. If $k=1$ then using the relations for $\langle a\,\underline{1}\,a\rangle $ and $\langle b\,\underline{1}\,b\rangle $ it is clear that any such path is a polynomial in $\langle aba\rangle$, hence is a polynomial in $u$. If $k>1$ then we can use the relation for $\langle \underline{k-1}\,\underline{k}\,\underline{k-1}\rangle$ to express the path as a linear combination of paths whose highest vertex number is $k-1$. The claim therefore follows by induction on $k$. Since $\mathcal{O}^{\boldsymbol{\lambda}}$ is isomorphic as a filtered algebra to some $D(\tilde{q})$ (by Theorem \ref{theo1}) it is not commutative and it is generated by $F_{2n-4}$. Now, the shortest paths starting and ending at $a$ and visiting $c$ or $d$ both have length $2n-4$; they are $\langle aca\rangle$ and $\langle ada \rangle$. Lemma \ref{recur1} tells us that $\langle aca \rangle + \langle ada \rangle$ is a polynomial in $u$. Thus $F_{2n-4}$ is spanned by $v$ and polynomials in $u$. Therefore $v$ cannot possibly commute with $u$ for otherwise $\mathcal{O}^{\boldsymbol{\lambda}}$ would be commutative.
\par
Next, we shall show that
\[ \left[ \langle aba \rangle,\langle aca \rangle \right] = \langle a(c-d) ba \rangle + \lambda_{ac}\langle aca \rangle - \lambda_{ad}\langle ada \rangle,\]
\[ \left[ \langle aba \rangle,\langle ada \rangle \right] = - \langle a(c-d)ba \rangle + \lambda_{ad}\langle ada \rangle - \lambda_{ac}\langle aca \rangle.\]
By symmetry it is sufficient to prove only the first of these. We shall repeatedly use Lemma \ref{zigzag} in what follows.
\begin{eqnarray*}
\langle abaca \rangle & = & \langle abdca \rangle - \lambda_{bc}\langle abca \rangle \\
& = & \langle abdba \rangle + \lambda_{ad} \langle abda \rangle - \lambda_{bc}\langle adca \rangle + \lambda_{bc}\lambda_{ac}\langle aca \rangle \\
& = & \langle acdba \rangle - \lambda_{ad}\langle adba \rangle + \lambda_{ad}\langle acda \rangle - \lambda_{ad}^2 \langle ada \rangle \\
& & -\lambda_{bc}\langle adba \rangle - \lambda_{bc}\lambda_{ad}\langle ada \rangle + \lambda_{bc}\lambda_{ac}\langle aca \rangle \\
& = & \langle acaba \rangle + \lambda_{bc}\langle acba \rangle - \lambda_{ad}\langle adba \rangle + \lambda_{ad}\langle acba \rangle \\
& & +\lambda_{ad}\lambda_{ac}\langle aca \rangle - \lambda_{ad}^2 \langle ada \rangle - \lambda_{bc}\langle adba \rangle - \lambda_{bc}\lambda_{ad}\langle ada \rangle +\lambda_{bc}\lambda_{ac}\langle aca \rangle \\
& = & \langle acaba \rangle + \langle acba \rangle - \langle adba \rangle + \lambda_{ac}\langle aca \rangle - \lambda_{ad}\langle ada \rangle,
\end{eqnarray*}
\noindent
using the fact that $\lambda_{ad}+\lambda_{bc} = 1$. Since $\left[ \langle aba \rangle,\langle aca \rangle \right] = \langle abaca \rangle - \langle acaba \rangle$ we are done. 
\par
Now we have
\begin{eqnarray*}
\left[ u,v \right] & = & \left[ \langle aba \rangle, \langle a(c-d)a \rangle \right] \\
& = & \left[ \langle aba \rangle, \langle aca \rangle \right] - \left[ \langle aba \rangle, \langle ada \rangle \right] \\
& = & 2\langle a(c-d)ba \rangle + 2\langle a (\lambda_{ac}c - \lambda_{ad}d)a\rangle \\
& = & 2\langle a(c-d)ba \rangle + 2\left\langle a \left(\frac{\lambda_{ac}-\lambda_{ad}}{2}(c+d)+\frac{\lambda_{ac}+\lambda_{ad}}{2}(c-d)\right)a\right\rangle \\
& = & 2\langle a(c-d)ba \rangle + 2\gamma\langle a(c+d)a\rangle + (1+2\mu_0)\langle a(c-d)a \rangle,
\end{eqnarray*}
\noindent
using $\frac{\lambda_{ac}-\lambda_{ad}}{2}=\gamma$, $\lambda_{ac}+\lambda_{ad}=1+2\mu_0$. Using the definitions of $v$, $w$ and $h(x)$ and the expression for $\langle a(c+d)a \rangle$ established in Lemma \ref{recur1} we get that
\begin{eqnarray*}\left[ u,v \right] & = & 2w - 2\mu_0 v - 2\left\{ ( \sqrt{u}-\mu_0 ) h(\sqrt{u})\right\} \\
& & +2\gamma\left\{ \prod_{i=0}^{n-2}( \sqrt{u}+\mu_{i}) \right\} + (1+2\mu_0)( v - \left\{ h(\sqrt{u})\right\}) \\
& = & 2w + v + \left\{ 2\gamma\prod_{i=0}^{n-2}\left( \sqrt{u}+\mu_i \right) - \left( 1+2\sqrt{u}\right)h\left(\sqrt{u}\right)\right\} \\
& = & 2w + v + \{\rho\} \\
& = & 2w + v. \hspace{5mm}\Box
\end{eqnarray*}

\begin{lemma}
\label{relat2}
In $\mathcal{O}^{\boldsymbol{\lambda}}$ we have $\left[u,w\right] = 2vu + w+ \rho$.
\end{lemma}

\noindent
Proof.
The proof is very similar to the proof of the previous lemma, although we still give all the details. First we shall prove that
\[ \left[ \langle aba \rangle,\langle acba \rangle \right] = \langle a(c-d) aba \rangle + \lambda_{bc}\langle acba \rangle - \lambda_{bd}\langle adba \rangle,\]
\[ \left[ \langle aba \rangle,\langle adba \rangle \right] = - \langle a(c-d)aba \rangle + \lambda_{bd}\langle adba \rangle - \lambda_{bc}\langle acba \rangle.\]
As before it is only necessary to prove the first of these. We have
\begin{eqnarray*}
\langle abacba \rangle & = & \langle abdcba \rangle - \lambda_{bc}\langle abcba \rangle \\
& = & \langle abdaba \rangle + \lambda_{bd} \langle abdba \rangle - \lambda_{bc}\langle adcba \rangle + \lambda_{bc}\lambda_{ac}\langle acba \rangle \\
& = & \langle acdaba \rangle - \lambda_{ad}\langle adaba \rangle + \lambda_{bd}\langle acdba \rangle - \lambda_{bd}\lambda_{ad} \langle adba \rangle \\
& & -\lambda_{bc}\langle adaba \rangle - \lambda_{bc}\lambda_{bd}\langle adba \rangle + \lambda_{bc}\lambda_{ac}\langle acba \rangle \\
& = & \langle acbaba \rangle + \lambda_{ac}\langle acaba \rangle - \lambda_{ad}\langle adaba \rangle + \lambda_{bd}\langle acaba \rangle \\
& & +\lambda_{bd}\lambda_{bc}\langle acba \rangle - \lambda_{bd}\lambda_{ad} \langle adba \rangle - \lambda_{bc}\langle adaba \rangle - \lambda_{bc}\lambda_{bd}\langle adba \rangle +\lambda_{bc}\lambda_{ac}\langle acba \rangle \\
& = & \langle acbaba \rangle + \langle acaba \rangle - \langle adaba \rangle + \lambda_{bc}\langle acba \rangle - \lambda_{bd}\langle adba \rangle,
\end{eqnarray*}
which is what we wanted. Next
\begin{eqnarray*}
\left[ u,w \right] & = & \left[ \langle aba \rangle, \langle a(c-d)ba \rangle \right] + \mu_0 \left[ u,v \right] \\
& = & \left[ \langle aba \rangle , \langle acba \rangle \right] - \left[ \langle aba \rangle , \langle adba \rangle \right] + \mu_0 (2w + v) \\
& = & 2\langle a(c-d)aba \rangle + 2 \langle a(\lambda_{bc}c - \lambda_{bd}d) ba \rangle + \mu_0 (2w + v) \\
& = & 2\langle a(c-d)a \rangle \langle aba \rangle + 2\gamma \langle a(c+d)ba \rangle + (1-2\mu_0)\langle a(c-d)ba \rangle + \mu_0 (2w+v),
\end{eqnarray*}
this time using $\frac{\lambda_{bc}-\lambda_{bd}}{2}=\gamma$, $\lambda_{bc}+\lambda_{bd}=1-2\mu_0.$
Using the definitions of $u$, $v$, $w$ and $h(x)$ and the expression for $\langle a(c+d)ba \rangle$ established in Lemma \ref{recur1} we get that
\begin{eqnarray*}
\left[ u,w \right] & = & 2\left( v - \left\{h(\sqrt{u})\right\}\right)(u - \mu_0^2) + 2\gamma\left\{ \left( \sqrt{u} - \mu_0 \right) \prod_{i=0}^{n-2}\left( \sqrt{u}+\mu_i \right) \right\} \\
& & +(1-2\mu_0)\left( w - \mu_0 v - \left\{ \left( \sqrt{u} - \mu_0\right) h\left( \sqrt{u}\right) \right\} \right) + \mu_0 (2w +v) \\
& = & 2vu + w + \left\{ \left( \sqrt{u} - \mu_0\right)\left( -\left( 1+2\sqrt{u}\right)h\left(\sqrt{u}\right)+2\gamma \prod_{i=0}^{n-2}\left( \sqrt{u} + \mu_i\right)\right)\right\} \\
& = & 2vu + w + \left\{ \left(\sqrt{u}-\mu_0\right)\rho\right\} \\
& = & 2vu + w + \rho.\hspace{5mm}\Box
\end{eqnarray*}

\par
Notice that Proposition \ref{mainprop} together with Lemmas \ref{relat1} and \ref{relat2} tell us that there is an isomorphism $\mathcal{O}^{\boldsymbol{\lambda}} \cong D(\tilde{q})$ such that the elements of $\mathcal{O}^{\boldsymbol{\lambda}}$ we have called $u$, $v$ and $w$ are identified with the corresponding elements of $D(\tilde{q})$. In particular we know that $ \mathcal{O}^{\boldsymbol{\lambda}} $ is free as a right $\mathbb{C} [u]$-module with free basis consisting of the monomials $v^i w^j$ where $j\leq 1$. It is therefore legitimate to use the notation $c(X)$ for the coefficient of $1$ (Definition \ref{const}). The problem is to prove that we can take $\tilde{q}$ to be $q$ as defined in Theorem \ref{theo2}.
 
\begin{defi}
For $0\leq k \leq n-3$ we define the following polynomials in $u$.
\[W_k:=c\left( \langle a(c-d) \,\underline{k} \,(c+d)a\rangle\right),\hspace{10mm} X_k:=c\left( \langle a(c-d) \,\underline{k} \,(c+d)ba\rangle\right),\]
\[Y_k:=c\left( \langle a(c-d) \,\underline{k} \,(c-d)a\rangle\right),\hspace{10mm}Z_k:=c\left( \langle a(c-d) \,\underline{k} \,(c-d)ba\rangle\right),\]
where $\underline{0}$ is interpreted as $(a+b)$. Also define
\[W_{n-2}:=c\left( \langle a(c-d)a \rangle \right),  \hspace{10mm} X_{n-2}:=c\left( \langle a(c-d)ba \rangle \right),\]
\[Y_{n-2}:=c\left( \langle a(c+d)a \rangle \right), \hspace{10mm}Z_{n-2}:=c\left( \langle a(c+d)ba \rangle \right).\]
\end{defi}

\begin{lemma}
\label{recur2}
For all $0\leq k \leq n-3$ we have
\[W_k=(\mu_{k+1}+\mu_0)W_{k+1}+X_{k+1},\hspace{10mm} X_k= (\mu_{k+1}-\mu_0)X_{k+1}+W_{k+1}(u-\mu_0^2),\]
\[Y_k=-2\gamma W_{k+1}-(1+\mu_0-\mu_{k+1})Y_{k+1}-Z_{k+1},\]
\[Z_{k}=-2\gamma X_{k+1}-(1-\mu_0 -\mu_{k+1})Z_{k+1}-Y_{k+1}(u-\mu_0^2).\]
\end{lemma}
\noindent
Proof. For brevity we shall define $H:=a(c-d)$. Suppose first that $0\leq k<n-3$. In what follows be shall apply Lemma \ref{kink} to $\langle \underline{k+1}\,\underline{k}\,\underline{n-3}\rangle$ and use the relations for $\langle c\,\underline{n-3}\,c\rangle$ and $\langle d\,\underline{n-3}\,d\rangle$. We have 
\begin{eqnarray*}
\langle H\underline{k}(c+d)a\rangle & = & \langle H\underline{k+1}\,\underline{k}\,\underline{n-3}\,(c+d)a\rangle \\
& = & \langle H \underline{k+1}\,(c+d)\,\underline{n-3}\,(c+d)a\rangle - \lambda_{k+1,n-3}\langle H \underline{k+1}\,\underline{n-3}\,(c+d)a\rangle \\
&=&\langle H \underline{k+1}(cda+dca)\rangle-\lambda_{k+1,n-3}\langle H\underline{k+1}(c+d)a\rangle\\
&&-\lambda_c \langle H\underline{k+1}ca\rangle - \lambda_d \langle H\underline{k+1}da\rangle.\\
\end{eqnarray*} 
By lemma \ref{zigzag} this equals
\[\langle H\underline{k+1}(c+d)ba\rangle + \lambda_{a,n-3}\langle H\underline{k+1}ca\rangle+ \lambda_{a,n-3}\langle H\underline{k+1}da\rangle-\lambda_{k+1,n-3}\langle H\underline{k+1}(c+d)a\rangle\]
\[=\langle H\underline{k+1}(c+d)ba\rangle + \lambda_{ak}\langle H\underline{k+1}(c+d)a\rangle.\]
By interchanging $a$ and $b$ in the above we also have
\[\langle H\underline{k}(c+d)b\rangle = \langle H\underline{k+1}(c+d)ab\rangle + \lambda_{bk}\langle H\underline{k+1}(c+d)b\rangle.\]
Now by multiplying on the right by $\langle ba\rangle$ we get
\[\langle H\underline{k}(c+d)ba\rangle = \langle H\underline{k+1}(c+d)a\rangle \langle aba\rangle + \lambda_{bk}\langle H\underline{k+1}(c+d)ba\rangle.\]
Since $\lambda_{ak}=\mu_{k+1}+\mu_0$, $\lambda_{bk}=\mu_{k+1}-\mu_0$ and $\langle aba\rangle=u-\mu_0^2$ we have therefore proved the relations for $W_k$ and $X_k$ for the case $0\leq k<n-3$. For the case $k=n-3$ essentially the same argument works but one has to set it out a little differently. We shall omit the details.
\par
Now we consider the relations for $Y_k$ and $Z_k$. Again we shall assume $0\leq k<n-3$ since one has to modify the argument slightly for the case $k=n-2$. As before we begin by applying Lemma \ref{kink}. We get
\begin{eqnarray*}
\langle H\underline{k}(c-d)a\rangle & = & \langle H \underline{k+1}\,(c+d)\,\underline{n-3}\,(c-d)a\rangle-\lambda_{k+1,n-3}\langle H\underline{k+1}(c-d)a\rangle\\
&=& \langle H \underline{k+1}(-cda+dca)\rangle - \lambda_{k+1,n-3}\langle H\underline{k+1}(c-d)a\rangle \\&&+ \langle H\underline{k+1}(-\lambda_c c +\lambda_d d)a\rangle\\
&=& \langle H\underline{k+1}(-cba+dba)\rangle - \lambda_{k+1,n-3}\langle H\underline{k+1}(c-d)a\rangle\\
&&+\langle H\underline{k+1}( (-\lambda_c -\lambda_{ac})c +(\lambda_d+\lambda_{ad})d)a\rangle.
\end{eqnarray*}
This is equal to
\[-\langle H\underline{k+1}(c-d)ba\rangle -\lambda_{k+1,n-3}\langle H\underline{k+1}(c-d)a\rangle\]
\[-\left(\frac{\lambda_c +\lambda_{ac}-\lambda_d - \lambda_{ad}}{2}\right)\langle H\underline{k+1}(c+d)a\rangle-\left(\frac{\lambda_c +\lambda_{ac}+\lambda_d + \lambda_{ad}}{2}\right)\langle H\underline{k+1}(c-d)a\rangle.\]
Next we use the fact that $\lambda_c +\lambda_{ac}-\lambda_d - \lambda_{ad}=4\gamma$ and that
\[\frac{\lambda_c +\lambda_{ac}+\lambda_d + \lambda_{ad}}{2}+\lambda_{k+1,n-3}=1-\lambda_{bk},\] 
to deduce that
\[\langle H\underline{k}(c-d)a\rangle = -2\gamma \langle H\underline{k+1}(c+d)a\rangle - (1-\lambda_{bk})\langle H\underline{k+1}(c-d)a\rangle - \langle H \underline{k+1}(c-d)ba\rangle.\]
The equations for both $Y_k$ and $Z_k$ follow from this, by repeating the argument at the end of the derivation of the equations for $W_k$ and $X_k$. $\Box$

\par
Before we solve these recurrence relations, it shall be convenient to introduce some notation. 
\begin{defi}
For $-1\leq k \leq n-2$ define $\displaystyle{p_k (x):=\prod_{i=k+1}^{n-2}(x + \mu_i)}$. Notice that these products do not involve $\mu_{n-1}$.
\end{defi}

\noindent
We now solve the recurrence relations.
\begin{lemma}
\label{recursolve}
For all $0\leq k \leq n-2$ we have
\begin{eqnarray*}
W_k & = & -\left\{ h( \sqrt{u}) p_k ( \sqrt{u} )\right\}, \\
X_k & = & -\left\{ ( \sqrt{u} -\mu_0 )h( \sqrt{u}) p_k ( \sqrt{u} )\right\}, \\
Y_k & = & 2\gamma\left\{ h( \sqrt{u}) \frac{p_k (\sqrt{u})-p_k ( -\sqrt{u}-1)}{1+2\sqrt{u}}\right\} + \left\{p_{-1}( \sqrt{u})p_k ( -\sqrt{u} -1)\right\},\\
Z_k & = & 2\gamma\left\{ ( \sqrt{u} -\mu_0) h( \sqrt{u}) \frac{p_k (\sqrt{u})-p_k ( -\sqrt{u}-1)}{1+2\sqrt{u}}\right\}\\
 && + \left\{(\sqrt{u}-\mu_0 )p_{-1}( \sqrt{u})p_k ( -\sqrt{u} -1)\right\}.
\end{eqnarray*}
\end{lemma}
\noindent
Proof. First note that $W_{n-2}=-\left\{ h(\sqrt{u}) \right\}$ and $X_{n-2}=-\left\{ ( \sqrt{u} - \mu_0 )h(\sqrt{u}) \right\}$ directly from the definitions of $v$ and $w$. Next $Y_{n-2}=\left\{p_{-1}(\sqrt{u})\right\}$ and $Z_{n-2}=\left\{(\sqrt{u}-\mu_0)p_{-1}(\sqrt{u})\right\}$ by Lemma \ref{recur1}. It is clear that substituting $k=n-2$ into the claimed expressions for $W_k$, $X_k$, $Y_k$ and $Z_k$ gives these results. Thus we only need to check that the expressions satisfy the recurrence relations of Lemma \ref{recur2}. Suppose $0\leq k < n-2$ and that the expressions for $W_{k+1}$, $X_{k+1}$, $Y_{k+1}$, $Z_{k+1}$ are correct. Then
\begin{eqnarray*}
X_k & = & (\mu_{k+1}-\mu_0)X_{k+1}+W_{k+1}(u-\mu_0^2) \\
& = & -\left\{ (\mu_{k+1}-\mu_0 )( \sqrt{u}-\mu_0)h( \sqrt{u})p_{k+1}( \sqrt{u})+(u-\mu_0^2)h( \sqrt{u})p_{k+1}( \sqrt{u})\right\}\\
& = & -\left\{ ( \sqrt{u} -\mu_0)h(\sqrt{u})(\mu_{k+1}-\mu_0 +\sqrt{u} +\mu_0)p_{k+1}( \sqrt{u})\right\} \\
& = & -\left\{ ( \sqrt{u} -\mu_0 )h( \sqrt{u}) p_k ( \sqrt{u} )\right\},
\end{eqnarray*}   
\noindent
as required. The proof for $W_k$ is similar so we omit the details. For $Y_k$ we get
\begin{eqnarray*}
Y_k & = & -2\gamma W_{k+1}-(1+\mu_0 - \mu_{k+1})Y_{k+1} -Z_{k+1} \\
& = & 2\gamma\left\{ \begin{matrix} h(\sqrt{u})p_{k+1}(\sqrt{u})-(1+\mu_0-\mu_{k+1})h( \sqrt{u})\frac{p_{k+1}( \sqrt{u})-p_{k+1}(-\sqrt{u}-1)}{1+2\sqrt{u}} \\ -\displaystyle{(\sqrt{u}-\mu_0)h( \sqrt{u})\frac{p_{k+1}( \sqrt{u})-p_{k+1}(-\sqrt{u}-1)}{1+2\sqrt{u}}} \end{matrix}\right\} \\
& & + \left\{ \begin{matrix} -(1+\mu_0 -\mu_{k+1})p_{-1}( \sqrt{u})p_{k+1}( -\sqrt{u}-1) \\ -(\sqrt{u}-\mu_0)p_{-1}(\sqrt{u})p_{k+1}(-\sqrt{u}-1) \end{matrix}\right\}.
\end{eqnarray*}
\noindent
In the first pair of braces, putting everything over a common denominator of $1+2\sqrt{u}$, the coefficient of $h(\sqrt{u})\frac{p_{k+1}( \sqrt{u})}{1+2\sqrt{u}}$ is $(1+2\sqrt{u}-1-\mu_0+\mu_{k+1}-\sqrt{u}+\mu_0)=(\sqrt{u}+\mu_{k+1})$, while the coefficient of $h(\sqrt{u})\frac{p_{k+1}( -\sqrt{u}-1)}{1+2\sqrt{u}}$ is $( 1+\mu_0 -\mu_{k+1}+\sqrt{u}-\mu_0)=-(-\sqrt{u}-1+\mu_{k+1})$, as required. The second pair of braces is easily seen to simplify to the required expression. The proof for $Z_k$ is done in exactly the same way. $\Box$ 
\par

We shall have no further use for Lemma \ref{recur2} and all we shall require from Lemma \ref{recursolve} are the expressions for $W_0$, $X_0$, $Y_0$ and $Z_0$. 

\begin{lemma} 
\label{magic} We have the following equality of polynomials.
\[p_{-1}(x)p_{-1}(-x-1)+h(x)h(-x-1)+\rho\frac{h(-x-1)-h(x)}{1+2x}=\frac{-4q(x)q(-x-1)+\rho^2}{(1+2x)^2}.\]
\end{lemma}

\noindent 
Proof. From the definition of $h$ we have 
\[(1+2x)h(x)=2\gamma p_{-1}(x)-\rho,\hspace{10mm}-(1+2x)h(-x-1)=2\gamma p_{-1}(-x-1)-\rho.\] Therefore we get 
\vspace{5mm}
\newline
$(1+2x)^2p_{-1}(x)p_{-1}(-x-1)+(1+2x)^2h(x)h(-x-1)+\rho (1+2x)(h(-x-1)-h(x))$
\begin{eqnarray*}
& = & (1+2x)^2p_{-1}(x)p_{-1}(-x-1)-(2\gamma p_{-1}(x)-\rho)(2\gamma p_{-1}(-x-1)-\rho) \\
&&+ \rho(-2\gamma p_{-1}(-x-1)+\rho - 2\gamma p_{-1}(x)+\rho) \\
& = & ((1+2x)^2-4\gamma^2)p_{-1}(x)p_{-1}(-x-1)+\rho^2 \\
& = & -4(x+\gamma +\tfrac{1}{2})(-x-1+\gamma +\tfrac{1}{2})p_{-1}(x)p_{-1}(-x-1)+\rho^2.
\end{eqnarray*}
The result follows because $\gamma+\tfrac{1}{2}=\mu_{n-1}$. $\Box$

\begin{lemma}
\label{wsquared}
In $\mathcal{O}^{\boldsymbol{\lambda}}$ we have 
\[c(w^2)-\mu_{0}c(wv)=\left\{( \sqrt{u}-\mu_0)\frac{-4q(\sqrt{u})q(-\sqrt{u}-1)+\rho^2}{(1+2\sqrt{u})^2}\right\}.\]
\end{lemma}

\noindent
Proof. By the definition of $w$ we have $\langle a(c-d)ba\rangle = w - \mu_{0} v - \left\{( \sqrt{u}-\mu_0 )h( \sqrt{u})\right\}$. Squaring both sides and taking the constant term we get
\begin{eqnarray*}
c\left( \langle a(c-d)ba(c-d)ba\rangle \right) & = & c(w^2) - \mu_0 c(wv)  - c\left( \left\{ ( \sqrt{u}-\mu_0 )h( \sqrt{u})\right\}w\right) \\
&& + \mu_0 c\left( \left\{ ( \sqrt{u}-\mu_0 )h( \sqrt{u})\right\}v\right) + \left\{ ( \sqrt{u}-\mu_0 )h( \sqrt{u})\right\}^2.
\end{eqnarray*}
On the other hand, we can use Lemma \ref{zigzag} repeatedly to deduce that
\vspace{5mm}
\newline
$\langle a(c-d)ba(c-d)ba \rangle$
\begin{eqnarray*} 
& = & \langle a(c-d)(bacba - badba)\rangle \\
& = & \langle a(c-d)(bdcba - \lambda_{bc} bcba - bcdba +\lambda_{bd} bdba )\rangle \\
& = & \langle a(c-d)(bdaba + 2\lambda_{bd} bdba  - bcaba - 2\lambda_{bc} bcba )\rangle \\
& = & -\langle a(c-d)b((c-d)aba + (\lambda_{bc} - \lambda_{bd})(c+d)ba + (\lambda_{bc} + \lambda_{bd})(c-d)ba)\rangle \\
& = & \langle a(c-d)a \rangle\langle a(c-d)aba + 2\gamma a(c+d)ba + (1-2\mu_0)a(c-d)ba\rangle \\
&& -\langle a(c-d)(a+b)( (c-d)aba + 2\gamma (c+d)ba + (1-2\mu_0)(c-d)ba)\rangle.
\end{eqnarray*}
Taking the constant term and using the definitions of $u$, $v$ and $w$ and the expression for $\langle a(c+d)ba\rangle$ derived in Lemma \ref{recur1} we get that $c\left( \langle a(c-d)ba(c-d)ba\rangle \right)$ is equal to 
\[c\left( ( v-\left\{h( \sqrt{u})\right\}) \left( \begin{matrix}
( v-\left\{h( \sqrt{u})\right\})(u-\mu_0^2) + 2\gamma\left\{ ( \sqrt{u}-\mu_0 )p_{-1}( \sqrt{u})\right\}\\+(1-2\mu_0 )( w - \mu_0 v - \left\{( \sqrt{u}-\mu_0 )h(\sqrt{u})\right\})\end{matrix}\right)\right)+M,\]
where $M := -Y_0 (u-\mu_0^2) -2\gamma X_0 - (1-2\mu_0)Z_0$. $c\left( \langle a(c-d)ba(c-d)ba\rangle \right)$ is then equal to
\[\left\{ h(\sqrt{u})\right\}^2 (u-\mu_0^2) - 2\gamma\left\{h(\sqrt{u})\right\}\left\{(\sqrt{u}-\mu_0)p_{-1}(\sqrt{u})\right\}\]
\[+ (1-2\mu_0)\left\{h(\sqrt{u})\right\}\left\{(\sqrt{u}-\mu_0)h( \sqrt{u})\right\}-c\left(\left\{h(\sqrt{u})\right\}v\right)(u-\mu_0^2)-(1-2\mu_0)c\left( \left\{h(\sqrt{u})\right\}w\right)\]
\[+\mu_0(1-2\mu_0)c\left( \left\{h(\sqrt{u})\right\}v\right)+M.\]
\noindent
Comparing our two expressions for $c\left( \langle a(c-d)ba(c-d)ba\rangle \right)$ we deduce that $c(w^2)-\mu_0 c(wv)$ is equal to $K + L + M$ where
\begin{eqnarray*}
K & := & -\left\{(\sqrt{u}-\mu_0)h( \sqrt{u})\right\}^2+\left\{ h(\sqrt{u})\right\}^2 (u-\mu_0^2) \\
& & - 2\gamma \left\{h(\sqrt{u})\right\}\left\{(\sqrt{u}-\mu_0)p_{-1}(\sqrt{u})\right\} +(1-2\mu_0)\left\{h(\sqrt{u})\right\}\left\{(\sqrt{u}-\mu_0)h(\sqrt{u})\right\},\\
L & := & -c\left( \left\{h(\sqrt{u})\right\}v\right)(u-\mu_0^2)-\mu_0 c\left(\left\{(\sqrt{u}-1+\mu_0)h(\sqrt{u})\right\}v\right) \\
&&+ c\left(\left\{(\sqrt{u}-1+\mu_0)h(\sqrt{u})\right\}w\right).
\end{eqnarray*}
\noindent
We shall now consider $K$, $L$ and $M$ in turn. For $K$, define $R:=\left\{h( \sqrt{u})\right\}$ and $S:=\left\{\sqrt{u}h( \sqrt{u})\right\}$. It is clear from the definition of $h(x)$ that $2\gamma p_{-1}(x)=(1+2x)h(x)+\rho$. Therefore 
\begin{eqnarray*}
-2\gamma\left\{(\sqrt{u}-\mu_0)p_{-1}(\sqrt{u})\right\} & = &-\left\{(\sqrt{u}-\mu_0)(1+2\sqrt{u})h(\sqrt{u})\right\}-\left\{(\sqrt{u}-\mu_0)\rho\right\} \\
&=&-((2u-\mu_0)R+(1-2\mu_0)S)-\rho.
\end{eqnarray*}
\noindent
We get 
\begin{eqnarray*}
K & = & -(S-\mu_0R)^2 + R^2 (u-\mu_0^2)-R((2u-\mu_0)R+(1-2\mu_0)S)-\rho R \\
&&+ (1-2\mu_0)R(S-\mu_0 R) \\
& = & -S^2 -uR^2 + 2\mu_0 RS - \rho R,
\end{eqnarray*}
which by Lemma \ref{bridge} is equal to $-\left\{ (\sqrt{u}-\mu_0)h(\sqrt{u})^2\right\}-\rho\left\{h(\sqrt{u})\right\}.$
\par
Now we consider $L$. By Lemma \ref{commutation} we can write $L$ in the form $\frac{\rho}{1-4u}\left\{\cdots\right\}$ where the contents of the braces is a $\mathbb{C} [\sqrt{u}]$ linear combination of $h( \sqrt{u} +1)$, $h( -\sqrt{u} -1)$ and $h(\sqrt{u})$. The fact that $1- 4u$ is not invertible in $\mathcal{O}^{\boldsymbol{\lambda}}$ is not an issue. One should recall that the expressions occurring in Lemma \ref{commutation} were in fact polynomials in $u$ and understand that here $u$ is being treated formally. Inside the braces, the coefficient of $h( \sqrt{u} +1)$ is
\[\frac{(1-2\sqrt{u})\left( -(u-\mu_0^2 ) - \mu_0( \sqrt{u}+\mu_0) +\sqrt{u}(\sqrt{u}+\mu_0)\right)}{2(\sqrt{u}+1)}=0.\]
The coefficient of $h( -\sqrt{u} -1)$ is
\[\frac{(1-2\sqrt{u})\bigg( (u-\mu_0^2 ) + \mu_0( -\sqrt{u}-2+\mu_0) -\sqrt{u}(-\sqrt{u}-2+\mu_0)\bigg)}{2(\sqrt{u}+1)}=(1-2\sqrt{u})(\sqrt{u}-\mu_0).\]
Finally, the coefficient of $h(\sqrt{u})$ is $-(2u-2\mu_0^2+2\mu_0(\sqrt{u}-1+\mu_0)+\sqrt{u}-1+\mu_0)=(1-2\sqrt{u})(\sqrt{u}+1+\mu_0)=(1-2\sqrt{u})(-(\sqrt{u}-\mu_0)+(1+2\sqrt{u}))$. We deduce that
\[L=\rho\left\{(\sqrt{u}-\mu_0)\frac{h(-\sqrt{u}-1)-h(\sqrt{u})}{1+2\sqrt{u}}+h(\sqrt{u})\right\}.\]
Finally we need to consider $M$. By Lemma \ref{recursolve} $M$ is equal to
\[-2\gamma\left\{h(\sqrt{u})\frac{p_{0}( \sqrt{u})-p_0 (-\sqrt{u}-1)}{1+2\sqrt{u}}\right\}(u-\mu_0^2) -\left\{p_{-1}(\sqrt{u})p_{0}(-\sqrt{u}-1)\right\}(u-\mu_0^2)\]
\[+2\gamma \left\{(\sqrt{u}-\mu_0)h(\sqrt{u})p_0 (\sqrt{u})\right\}-2\gamma(1-2\mu_0)\left\{(\sqrt{u}-\mu_0)h(\sqrt{u})\frac{p_{0}( \sqrt{u})-p_0 (-\sqrt{u}-1)}{1+2\sqrt{u}}\right\}\]
\[-(1-2\mu_0)\left\{(\sqrt{u}-\mu_0)p_{-1}(\sqrt{u})p_0 (-\sqrt{u}-1)\right\}.\]
The two terms not involving $2\gamma$ combine to give \[-\left\{(\sqrt{u}-\mu_0)(\sqrt{u}-\mu_0+1)p_{-1}(\sqrt{u})p_0 (-\sqrt{u}-1)\right\} = \left\{(\sqrt{u}-\mu_0)p_{-1}(\sqrt{u})p_{-1}(-\sqrt{u}-1)\right\}.\]
Similarly, the three remaining terms combine to give
\[2\gamma\left\{(\sqrt{u}-\mu_0)h(\sqrt{u})\frac{p_{-1}( \sqrt{u})-p_{-1}( -\sqrt{u}-1)}{1+2\sqrt{u}}\right\}.\]
Using again the fact that $2\gamma p_{-1}(x)=(1+2x)h(x)+\rho$ this simplifies to 
\[\left\{(\sqrt{u}-\mu_0)\bigg((h\sqrt{u})^2+(h\sqrt{u})h(-\sqrt{u}-1)\bigg)\right\}.\]
We therefore deduce that $c(w^2)-\mu_0 c(wv)=K + L + M$ is equal to
\[\left\{( \sqrt{u}-\mu_0)\left(p_{-1}(\sqrt{u})p_{-1}(-\sqrt{u}-1)+h(\sqrt{u})h(-\sqrt{u}-1)+\rho\frac{h(-\sqrt{u}-1)-h(\sqrt{u})}{1+2\sqrt{u}}\right)\right\}.\]
The result therefore follows from Lemma \ref{magic}. $\Box$

We are now in a position to prove Theorem \ref{theo2}. Using Proposition \ref{mainprop} and Lemmas \ref{relat1}, \ref{relat2} we know that there is an isomorphism of $\mathcal{O}^{\boldsymbol{\lambda}}$ with some $D(\tilde{q})$ such that the elements we have called $u$, $v$, $w$ are identified with the corresponding elements of $D(\tilde{q})$. By considering the relations for $D(\tilde{q})$ we see that there exists a polynomial $p(x)$ satisfying $p(-x)=p(x-1)$ such that $c(w^2)=\left\{\sqrt{u}p(\sqrt{u})\right\}$ and $c(wv)=\left\{p(\sqrt{u})\right\}$. It follows that $c(w^2)-\mu_0 c(wv)=\left\{(\sqrt{u}-\mu_0)p(\sqrt{u})\right\}$. Since this property holds for the polynomial \[\frac{-4q(x)q(-x-1)+\rho^2}{(1+2x)^2},\] Lemma \ref{iso} tells us that $p(x)$ equals this polynomial. Therefore $u$, $v$, $w$ must satisfy the desired relations. $\Box$

\section{Reflections and Triality}

\begin{defi}
For any vertex $j$ of $Q$ we define a linear map $r_j :\mathbb{C}^I \rightarrow \mathbb{C}^I$ by 
\[ (r_j (\boldsymbol{\nu}))^i =\begin{cases} -\boldsymbol{\nu}^i & \quad i=j \\
\hfill \boldsymbol{\nu}^i+\boldsymbol{\nu}^j & \quad i-j \\  \boldsymbol{\nu}^i & \text{otherwise}, \end{cases}
\]
where $i-j$ means that $i$, $j$ are connected in the underlying graph of $Q$.  
\end{defi}

\begin{theo}
\label{theo3}
Let $n\geq 4$. Every deformation of the coordinate algebra of a Kleinian singularity of type $D_n$ that is not commutative is isomorphic as a filtered algebra to $\mathcal{O}^{\boldsymbol{\lambda}}$ for some $\boldsymbol{\lambda}$ satisfying $\boldsymbol{\lambda}\cdot\boldsymbol{\delta}=1$. If $\boldsymbol{\lambda}$ and $\boldsymbol{\tilde{\lambda}}$ are such that $\boldsymbol{\lambda}\cdot\boldsymbol{\delta}=\boldsymbol{\tilde{\lambda}}\cdot\boldsymbol{\delta}=1$ then $\mathcal{O}^{\boldsymbol{\lambda}}\cong\mathcal{O}^{\boldsymbol{\tilde{\lambda}}}$ as filtered algebras if and only if $\boldsymbol{\tilde{\lambda}}$ can be obtained from $\boldsymbol{\lambda}$ by applying a sequence of reflections $r_j$ $(j\neq a)$ and graph automorphisms fixing $a$.
\end{theo}
\noindent 
Proof. The first statement follows from Theorems \ref{theo1} and \ref{theo2} since if $q(x)=\prod_{i=0}^{n-1}(x+\mu_i)$ then $D(q)\cong \mathcal{O}^{\boldsymbol{\lambda}}$ where $\boldsymbol{\lambda}$ is defined to be
\[
\xymatrix{ (\mu_{0}+\mu_1) \ar@{-}[rd] & &  & & (\mu_{n-1}-\mu_{n-2}) \\ & (\mu_{2}-\mu_1)\ar@{-}[r] & (\mu_3-\mu_2)\ar@{.}[r] & (\mu_{n-2}-\mu_{n-3})\ar@{-}[ur]\ar@{-}[dr] & \\  (\mu_1 - \mu_0)\ar@{-}[ur] & &  & & (1-\mu_{n-1}-\mu_{n-2}).}
\]
To any $\boldsymbol{\lambda}\in\mathbb{C}^I$ we associate the ordered sequence $(\mu_0 , \mu_1, \cdots, \mu_{n-1})$, where the entries in this sequence are defined as in Section 6. It follows from Theorem \ref{theo1} that if $n\geq5$,  $\mathcal{O}^{\boldsymbol{\lambda}}\cong\mathcal{O}^{\boldsymbol{\tilde{\lambda}}}$ if and only if the sequence for $\boldsymbol{\tilde{\lambda}}$ can be obtained from the sequence for $\boldsymbol{\lambda}$ by applying a permutation of the entries and possibly replacing some of the $\mu_i$ by $1-\mu_i$. However, it is easy to see that if we apply $r_b$ to the vector above we get  
\[
\xymatrix{ (\mu_{0}+\mu_1) \ar@{-}[rd] & &  & & (\mu_{n-1}-\mu_{n-2}) \\ & (\mu_{2}-\mu_0)\ar@{-}[r] & (\mu_3-\mu_2)\ar@{.}[r] & (\mu_{n-2}-\mu_{n-3})\ar@{-}[ur]\ar@{-}[dr] & \\  (\mu_0 - \mu_1)\ar@{-}[ur] & &  & & (1-\mu_{n-1}-\mu_{n-2}).}
\]
\noindent
Notice that $r_b$ has simply switched the first two entries of the sequence. Similarly we see that $r_1$, ..., $r_{n-3}$, $r_c$ correspond to the remaining elementary transpositions. $r_d$ on the other hand corresponds to the transformation $(\mu_0, \cdots, \mu_{n-3}, \mu_{n-2},\mu_{n-1})\mapsto (\mu_0, \cdots, \mu_{n-3}, 1-\mu_{n-1},1-\mu_{n-2})$. Therefore using the reflections $r_i$ for $i\neq a$ we can permute the entries of the sequence anyhow we like and replace an \emph{even} number of the $\mu_i$ by $1-\mu_i$. However the graph automorphism switching $c$ and $d$ has the effect of simply replacing $\mu_{n-1}$ by $1-\mu_{n-1}$. For $n\geq 5$ the result follows.
\par
For $n=4$ the situation is complicated by the fact that there are more graph automorphisms. If we apply the graph automorphism switching $b$ and $c$ to the vector $\boldsymbol{\lambda}$ and work out the corresponding sequence we get
\[\left(\frac{\mu_0+\mu_1+\mu_2-\mu_3}{2}, \frac{\mu_0+\mu_1-\mu_2+\mu_3}{2},\frac{\mu_0-\mu_1+\mu_2+\mu_3}{2},\frac{-\mu_0+\mu_1+\mu_2+\mu_3}{2}\right).\]
This corresponds to an isomorphism $D(q')\rightarrow \mathcal{O}^{\boldsymbol{\lambda}}$ where \[q'(x)=\textstyle{\left(x+\frac{\mu_0+\mu_1+\mu_2-\mu_3}{2}\right)\left(x+\frac{\mu_0+\mu_1-\mu_2+\mu_3}{2}\right)\left(x+\frac{\mu_0-\mu_1+\mu_2+\mu_3}{2}\right)\left(x+\frac{-\mu_0+\mu_1+\mu_2+\mu_3}{2}\right)}.\] 
By considering Remark \ref{orientation}, it is not too difficult to see that the image of the first standard generator of $D(q')$ is $u':=-\langle aca\rangle + \textstyle{\left( \frac{\mu_0 +\mu_1 +\mu_2 - \mu_3}{2}\right)^2}$. Similarly the graph automorphism switching $b$ and $d$ corresponds to an isomorphism $D(q'')\rightarrow \mathcal{O}^{\boldsymbol{\lambda}}$ where
\[q''(x)=\textstyle{\left(x+\frac{\mu_0+\mu_1+\mu_2+\mu_3 -1}{2}\right)\!\left(x+\frac{\mu_0+\mu_1-\mu_2-\mu_3 +1}{2}\right)\!\left(x+\frac{\mu_0-\mu_1+\mu_2-\mu_3 +1}{2}\right)\!\left(x+\frac{\mu_0 -\mu_1 -\mu_2+\mu_3 +1}{2}\right)}.\]
The image of the first standard generator of $D(q'')$ is $u'':=-\langle ada\rangle + \textstyle{\left( \frac{\mu_0 +\mu_1 +\mu_2 + \mu_3 -1}{2}\right)^2}$. However from the proofs of Proposition \ref{mainprop} and Theorem \ref{theo1} we know that there are precisely three pairs $(\tilde{u},\tilde{v})$ of elements of $\mathcal{O}^{\boldsymbol{\lambda}}$ that correspond to an isomorphism of $\mathcal{O}^{\boldsymbol{\lambda}}$ with some $D(\tilde{q})$ (up to multiplying $\tilde{v}$ by a nonzero scalar). Since $\mathcal{O}^{\boldsymbol{\lambda}}$ is a deformation of a Kleinian singularity of type $D_4$, the dimension of the space $F_4$ is known to be $3$. In fact it is easy to see that it is spanned by $1$, $\langle aba \rangle$, $\langle aca \rangle$ and $\langle ada \rangle$ and that these satisfy the unique relation $\langle aca\rangle + \langle ada \rangle = \langle aba \rangle + \lambda_a(\lambda_a + \lambda_1)1$. Hence the elements $u$, $u'$ and $u''$ are all different. It follows that $D(q)\cong D(\tilde{q})$ if and only if $\tilde{q}(x)\tilde{q}(-x-1)$ is a nonzero scalar multiple of one of $q(x)q(-x-1)$, $q'(x)q'(-x-1)$, $q''(x)q''(-x-1)$. Thus $\mathcal{O}^{\boldsymbol{\tilde{\lambda}}}\cong \mathcal{O}^{\boldsymbol{\lambda}}$ if and only if the sequence associated to $\boldsymbol{\tilde{\lambda}}$ is obtained from one of the three sequences $(\mu_0 ,\mu_1, \mu_2, \mu_3)$, $\left(\frac{\mu_0+\mu_1+\mu_2-\mu_3}{2}, \frac{\mu_0+\mu_1-\mu_2+\mu_3}{2},\frac{\mu_0-\mu_1+\mu_2+\mu_3}{2},\frac{-\mu_0+\mu_1+\mu_2+\mu_3}{2}\right)$ or\newline $\left(\frac{\mu_0\!+\!\mu_1\!+\!\mu_2\!+\!\mu_3\!-\!1}{2}, \frac{\mu_0\!+\!\mu_1\!-\!\mu_2\!-\!\mu_3\!+\!1}{2},\frac{\mu_0\!-\!\mu_1\!+\!\mu_2\!-\!\mu_3\!+\!1}{2},\frac{\mu_0\!-\!\mu_1\!-\!\mu_2\!+\!\mu_3\!+\!1}{2}\right)$ by permuting the entries and possibly transforming one or more of the entries by $\mu\mapsto 1-\mu$. The result follows. $\Box$

Email: paulsboddington@yahoo.co.uk

\begin{thebibliography}{99}
\bibitem{a} Arakawa, T.: Representation Theory of W-algebras, arXiv:math.QA/0506056
\bibitem{bav1} Bavula, V.: Finite-dimensionality of {${\rm Ext}\sp n$} and {${\rm Tor}\sb n$} of simple modules over a class of algebras, Funktsional. Anal. i Prilozhen, 25(3), 80-82, 1991.
\bibitem{bav2} Bavula, V.: Generalized {W}eyl algebras and their representations, Algebra i Analiz, 4(1): 75-97, 1992.
\bibitem{bav} Bavula, V. and Jordan, D. A.: Isomorphism Problems and Groups of Automorphisms for Generalized Weyl Algebras, Tran. Amer. Math. Soc {\bf353}, no. 2 (2000) 769-794
\bibitem{boddington} Boddington, P.: No-Cycle Algebras and Representation Theory, PhD Thesis, University of Warwick (2004)
\bibitem{brieskorn} Brieskorn, E.: Singular Elements of Semisimple Algebraic Groups, Actes Congrs Intern. Math {\bf2} (1970) 279-284
\bibitem{CBH} Crawley-Boevey, W. and Holland, M.: Noncommutative Deformations of Kleinian Singularities, Duke Math J, {\bf92}, no. 3 (1998) 605-635
\bibitem{boer} de Boer, J. and Tjin, T.: Quantizations and Representation Theory of Finite W-algebras, Commun. Math. Phys. {\bf158} (1993) 485-516
\bibitem{SK} De Sole, A. and Kac, V. G.: Finite vs Affine W-Alegbras, arXiv:math-ph/0511055
\bibitem{gordonrumynin} Gordon, I. and Rumynin, D.: Subregular Representations of $\mathfrak{sl}_n$ and Simple Singularities of Type $A_{n-1}$, Compositio Math. {\bf138} (2003) 337-360
\bibitem{hodges} Hodges, T. J.: Non-commutative Deformations of Type-A Kleinian Singularities, J. Alg. {\bf161} (1993) 271-290
\bibitem{levy} Levy, P.: Isomorphism Problems of Noncommutative Deformations of Type $D$ Kleinian Singularities, arXiv:math.RA/0610490 (2006)
\bibitem{mckay} McKay, J.: `Graphs, Singularities, and Finite Groups' in The Santa Cruz Conference on Finite Groups (Univ. California, Santa Cruz, Calif., 1970) Proc. Sympos. Pure Math. {\bf 37}, Amer. Math. Soc. Providence (1980) 183-186
\bibitem{premet} Premet, A.: Special Transverse Slices and Their Enveloping Algebras, Adv. Math. (2002) 1-55 
\bibitem{slodowy} Slodowy, P.: Simple Singularities and Simple Algebraic Groups, Lecture Notes in Mathematics, Vol. 815, Springer, Berlin/Heidelberg/New York (1980)
\bibitem{smith} Smith, S.P.: A Class of Algebras Similar to the Enveloping Algebra of $\mathfrak{sl}(2)$, Trans. Amer. Math. Soc., {\bf322} (1990) 285-314
\end{thebibliography}
\end{document}